\documentclass[11pt, a4paper]{amsart}

\usepackage{amsmath, amssymb, amsthm}
\usepackage[english]{babel}
\usepackage[all]{xy}
\usepackage{graphicx}
\usepackage{tikz}
\usetikzlibrary{matrix}
\usepackage{verbatim}
\usepackage{enumerate}
\usepackage{url}
\usepackage{subfig}
\usepackage{bm}
\usepackage{xcolor, overpic}

\newtheorem{thm}{Theorem}[section]

\theoremstyle{definition}

\theoremstyle{remark}
\newtheorem{rem}[thm]{Remark}

\newcommand{\R}{\mathbb{R}}

\newcommand{\note}[1]{}

\usepackage{xcolor}
\definecolor{colCB}{rgb}{0,0,.7}
\definecolor{colTB}{rgb}{.7,0,0}
\definecolor{colOO}{rgb}{0,.7,0}

\newcommand{\T}{\mathbb T}
\newcommand{\Z}{\mathbb Z}
\newcommand{\N}{\mathbb N}
\renewcommand{\S}{\mathbb S}

\renewcommand{\d}{\mathrm d}
\newcommand{\abs}[1]{\lvert #1 \rvert}

\newcommand{\fr}[2]{\frac{\displaystyle #1}{\displaystyle #2}}
\newcommand{\df}[2]{\frac{\displaystyle d#1}{\displaystyle d#2}}
\newcommand{\pf}[2]{\frac{\displaystyle \partial #1}{\displaystyle \partial #2}}

\newcommand{\Real}{\mathrm{Re}\:}
\newcommand{\Imag}{\mathrm{Im}\:}

\begin{document}



\title[Hopf Bifurcations of Twisted States]{
Hopf Bifurcations of Twisted States in Phase Oscillators Rings with Nonpairwise Higher-Order Interactions}
\author[Christian Bick, Tobias B\"ohle,  Oleh E.~Omel'chenko]{Christian Bick${}^1$, Tobias B\"ohle${}^2$, and Oleh E.~Omel'chenko${}^3$}%

\address{${}^1$Department of Mathematics, Vrije Universiteit Amsterdam, DE Boelelaan 1111, Amsterdam, The Netherlands;
Institute for Advanced Study, Technical University of Munich, Lichtenbergstr 2, 85748 Garching, Germany; 
Department of Mathematics, University of Exeter, Exeter EX4 4QF, United Kingdom; 
and Mathematical Institute,
University of Oxford, Oxford OX2 6GG, United Kingdom.}
\address{$^{2}$Technical University of Munich, School of Computation Information and Technology, Department of Mathematics, Boltzmannstr 3, 85748 Garching, Germany;
Institute for Advanced Study, Technical University of Munich, Lichtenbergstr 2, 85748 Garching, Germany; 
and Department of Mathematics, Vrije Universiteit Amsterdam, DE Boelelaan 1111, Amsterdam, The Netherlands
(tobias.boehle@tum.de).}
\address{${}^3$Institute of Physics and Astronomy, University of Potsdam, Karl-Liebknecht-Str. 24/25, 14476 Potsdam, Germany}
\date{\today}

\maketitle

\begin{abstract}
Synchronization is an essential collective phenomenon in networks of interacting oscillators. 
Twisted states are rotating wave solutions in ring networks where the oscillator phases wrap around the circle in a linear fashion. 
Here, we analyze Hopf bifurcations of twisted states in ring networks of phase oscillators with nonpairwise higher-order interactions. 
Hopf bifurcations give rise to quasiperiodic solutions that move along the oscillator ring at nontrivial speed.
Because of the higher-order interactions, these emerging solutions may be stable.
Using the Ott--Antonsen approach, we continue the emergent solution branches which approach anti-phase type solutions (where oscillators form two clusters whose phase is~$\pi$ apart) as well as twisted states with a different winding number.
\end{abstract}

\section{Introduction}
\begin{figure}
\begin{center}
\includegraphics[width=0.32\textwidth]{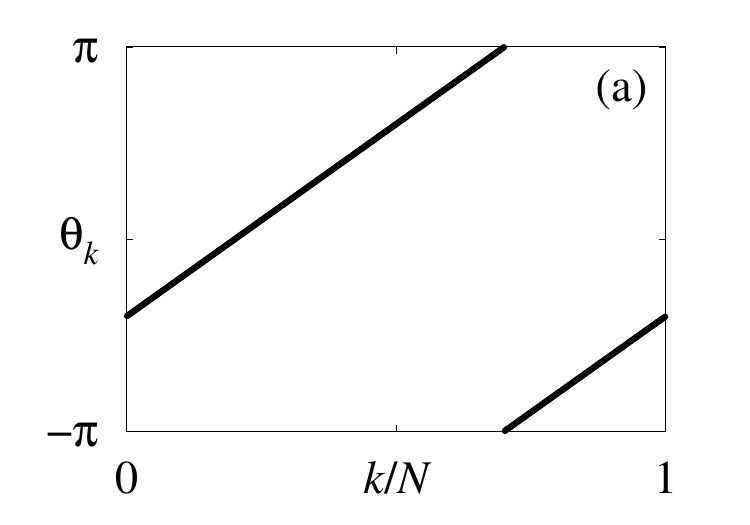}
\includegraphics[width=0.32\textwidth]{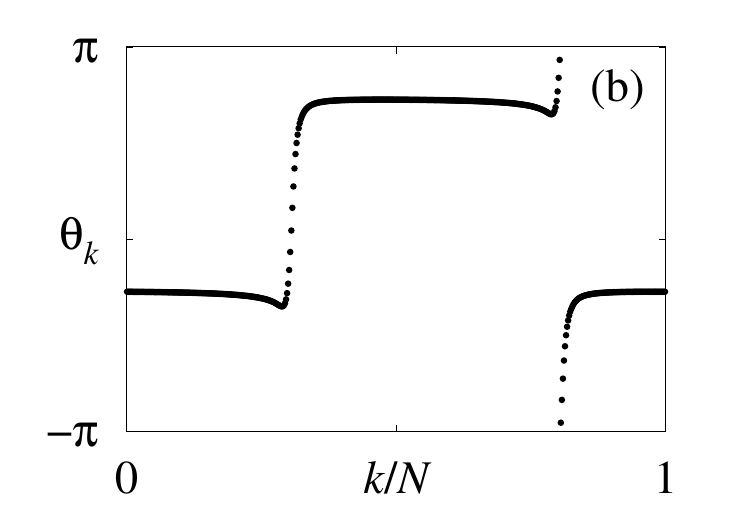}
\includegraphics[width=0.32\textwidth]{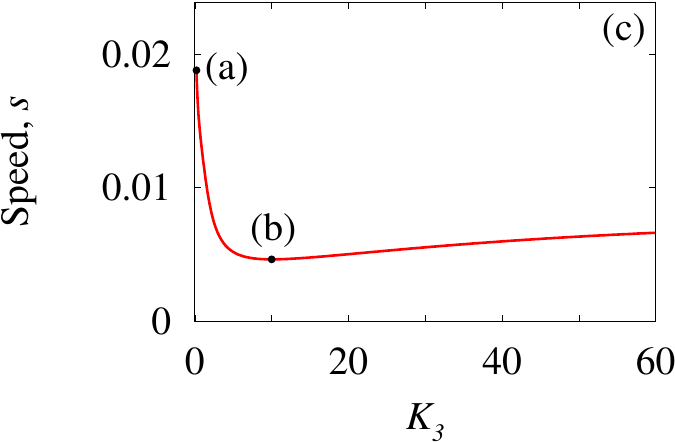}
\end{center}
\caption{
\label{fig:Teaser}
Twisted states in oscillator rings~\eqref{eq:finite_system} of $N = 512$ oscillators with distance-dependent coupling~\eqref{Formula:G}.
Panel~(a) shows a $1$-twisted (splay) state for $K_3 \approx 0.22$ and Panel~(b) a traveling nonuniformly twisted wave for $K_3 = 10$ for $\omega = 0$, $A = 0.9$, $B = 0.1$, $K_2 = 1$, $\alpha_2 = \pi/2 - 0.1$ and $\alpha_3 = 0$.
Panel~(c) shows theoretically predicted dependence of the drift speed $s$ of traveling waves versus the parameter $K_3$ for $B = 0.1$.
The dots indicate values corresponding to phase snapshots in Panels~(a) and~(b).
}
\end{figure}

\noindent
Coupled phase oscillators on a network provide essential models to understand synchronization phenomena~\cite{Rodrigues2016}. 
Apart from global phase synchrony, where the phase of all oscillators coincides, twisted states in nonlocally coupled networks have attracted attention~\cite{Wiley2006}. 
For example, consider $N$~Kuramoto oscillators on a ring network, whose phase $\theta_k\in\T := \R/2\pi\Z$ evolves according to
\begin{align}\label{eq:finite_system_kur}
\begin{split}
	\dot \theta_k := \frac{\d}{\d t}\theta_k &= \omega + \frac{1}{N} \sum_{j=1}^N G\left( \frac{k-j}{N}\right) \sin(\theta_j-\theta_k),
\end{split}
\end{align}
where $\omega\in\R$ is the intrinsic oscillator frequency and $G:\R\to\R$ is a smooth one-periodic function that determines the strength of the oscillator interaction depending on their relative position on the ring. 
For $q\in\N$, the \emph{$q$-twisted states} 
\begin{align}\label{eq:SplayFinite}
	\Theta_k^q(t) := 2 \pi q \frac{k}{N} + \Omega t + \beta,
\end{align}
with $\beta\in\T$, $\Omega\in\R$, $k=1,\dots,N$ are periodic solutions of~\eqref{eq:finite_system_kur}; see Figure~\ref{fig:Teaser}(a) for an example.
These solutions are also known as rotating wave solutions or, for $q=1$, as splay phase configurations. 
While the stability of such solutions has been analyzed explicitly~\cite{Wiley2006,Girnyk2012}, the specific form of Kuramoto phase coupling (a single harmonic without phase shift) imposes a gradient structure, which prevents the emergence of bifurcations to periodic solutions.

Richer dynamics and bifurcation behavior are possible for more generic phase interactions that one would expect from phase reductions~\cite{Nakao2015}. 
Indeed, phase reductions also give rise to nonpairwise phase interaction terms~\cite{Ashwin2015a,Leon2019a,Bick2023}, which can give rise to Hopf bifurcations of splay phase solutions and more complicated dynamical behavior~\cite{Bick2016b}. 
A generalization of the Kuramoto model~\eqref{eq:finite_system_kur} for ring-like networks to include nonpairwise interaction terms is
\begin{align}\label{eq:finite_system}
\begin{split}
	\dot \theta_k &= \omega + \frac{K_2}{N} \sum_{j=1}^N G\left( \frac{k-j}{N}\right) \sin(\theta_j-\theta_k+\alpha_2)\\
	&\quad\  + \frac{K_3}{N^2}\sum_{k,l=1}^N G\left(\frac{k-j}{N}\right) G\left( \frac{k-l}{N}\right) \sin(2 \theta_j - \theta_l-\theta_k + \alpha_3),
\end{split}
\end{align}
where~$K_2\in\R$ determines the strength of the pairwise interactions with a phase shift~$\alpha_2\in\T$ and $K_3\in \R$ determines the strength of nonlinear interactions between three phase variables. Such phase oscillator networks with ``higher-order'' nonpairwise interactions have their intrinsic interest~\cite{Battiston2020,Bick2021} as dynamical systems on (weighted) hypergraphs where the \emph{pairwise phase interaction function} $\sin(\theta_j-\theta_k+\alpha_2)$ corresponds to interactions along edges and the \emph{triplet phase interaction function} $\sin(2 \theta_j - \theta_l-\theta_k + \alpha_3)$ corresponds to interactions along hyperedges.

The typical setup to analyze $q$-twisted states is the continuum limit of $N\to\infty$ oscillators; see already~\cite{Wiley2006}. 
In the continuum limit the phase of oscillator $x\in\S := [0,1]/(1\sim 0)$---the unit interval with end points identified---at time~$t$ is given by a phase $\theta(x,t) = \theta_x(t)$ that, for the ring network~\eqref{eq:finite_system}, evolves according to
\begin{equation}\label{eq:main_system_w}
\begin{split}
        \dot \theta_x &= \omega + K_2 \int_\S G(x-y) \sin(\theta_y-\theta_x+\alpha_2)\ \d y\\
	& \quad\ + K_3 \int_\S \int_\S G(x-y)G(x-z) \sin(2 \theta_y - \theta_z-\theta_x + \alpha_3)\ \d z \d y.
\end{split}    
\end{equation}
Now $q$-twisted states 
\begin{align}\label{eq:SplayContinuum}
	\Theta^q(x,t) := 2 \pi q x + \Omega t + \beta
\end{align}
are periodic solutions with a smooth phase profile (in~$x$). Their stability has been analyzed for networks with Kuramoto coupling ($\alpha_2=0$) without higher-order interactions~\cite{Wiley2006,Medvedev2015} and more recently also with higher-order interaction~\cite{Bick2023a}.

Here we analyze Hopf bifurcations of $q$-twisted states in phase oscillator networks~\eqref{eq:main_system_w} with higher-order interactions. 
To identify Hopf bifurcation points, we linearize the system at $q$-twisted states and find conditions for a complex conjugate pair of eigenvalues to cross the imaginary axis; we focus primarily on $1$-twisted states. 
Computing higher-order derivatives allows to determine whether the Hopf bifurcation is subcritical or supercritical and estimate amplitude and period of the bifurcating solution. 
The emergent periodic solutions are quasiperiodic solutions with nontrivial rotation along the spatial domain~$\S$---see Figure~\ref{fig:Teaser}(b) for an example of the corresponding solution in a finite network. Because of the higher-order interactions, these emergent solutions can be stable. To continue the solutions further from the bifurcation point, we consider the dynamics on the Ott--Antonsen manifold~\cite{Ott2008,Bick2018c}.
Adapting recent approaches to continue periodic solutions on the Ott-Antonsen manifold~\cite{Omelchenko2022,Omelchenko2023}, we compute bifurcation branches for varying parameters; cf.~Figure~\ref{fig:Teaser}(c) for an example of branch of periodic solutions bifurcating from a $1$-twisted state.
Further bifurcations involve (traveling and stationary) antiphase solutions and some branches appear to come close to $-1$-twisted solution.

This paper is organized as follows. In Section~\ref{sec:Prelim} we specify the dynamical equations as well as relevant parameters and discuss them in the context of recent related work. In Section~\ref{sec:bifurcation}, we linearize the equations and give the eigenvalues that determine the bifurcations as well as the higher-order derivatives that determine the bifurcation type.
We consider the system on the Ott--Antonsen manifold in Section~\ref{sec:OttAnt} and outline the continuation technique. In Section~\ref{sec:Continuation} we describe the bifurcations that arise and conclude in Section~\ref{sec:Discussion} with some remarks.

\section{Twisted States in Ring Networks with Higher-Order Interactions}
\label{sec:Prelim}
\noindent
In the following we analyze $q$-twisted states in the continuum limit equation~\eqref{eq:main_system_w}. Without loss of generality we will henceforth set $\omega = 0$ by exploiting the phase-shift symmetry $\theta_x \mapsto \theta_x+\beta$, $\beta\in\T$, that shifts all oscillator phases by a constant phase angle~$\beta$. Thus, the dynamical equations are
\begin{subequations}\label{eq:main_system}
\begin{align}
    \dot \theta_x &= K_2 \int_\S G(x-y) \sin(\theta_y-\theta_x+\alpha_2)\ \d y\\
	\label{eq:main_system_triplet}
	& \quad\ + K_3 \int_\S \int_\S G(x-y)G(x-z) \sin(2 \theta_y - \theta_z-\theta_x + \alpha_3)\ \d z \d y,
\end{align}
\end{subequations}
and the $q$-twisted states~\eqref{eq:SplayContinuum} are actually equilibria relative to the phase-shift symmetry action. The $0$-twisted state corresponds to full phase synchrony in the system and the $\pm1$-twisted states to the classical (anti)splay phase configuration.

The network has a ring structure (with orientation) since the system has a spatial rotational symmetry, where $s\in\S$ acts on~$\S$ by $s:x\mapsto x+s$. Specifically, we consider a network with a coupling kernel function
\begin{equation}
G(x) = 1 + A \cos(2 \pi x) + B \sin(2 \pi x).
\label{Formula:G}
\end{equation}
with only the first nontrivial harmonic being present; see also~\cite{Abrams2004,Omelchenko2018,Omelchenko2020}. A nonzero parameter~$B$ breaks the reflectional symmetry of the ring network: For $B=0$, the reflection $x\mapsto -x$ is also a symmetry of the system.

The higher-order triplet interactions are determined by the triplet coupling kernel
\begin{equation}\label{eq:TripletKernel}
W(x,y,z)=G(x-y)G(x-z)    
\end{equation}
that is a product of the same coupling kernel function~$G$ that determines the pairwise interactions. Such a product structure is natural if the phase equations are obtained through a phase reduction~\cite{Bick2023}. However, it is distinct from other triplet interaction kernels considered in the literature in the context of ring networks. 
Specifically, in~\cite{Bick2023a} the authors considered a generalized top-hat coupling kernel $W(x,y,z) = W_r(z+y-2x)$ where $W_r(v) =  1$ if $\min(\abs{v}, 1-\abs{v})\le r$ and $W_r(v) = 0$ otherwise.
This interaction function lacks a product structure but generalizes coupling with a finite coupling range considered in the context of twisted states on rings with pairwise coupling~\cite{Girnyk2012}.
By contrast, with the triplet coupling kernel~\eqref{eq:TripletKernel} and the kernel function~\eqref{Formula:G}, the interactions in~\eqref{eq:main_system} only consists of finitely many Fourier modes; this facilitates analytic computations, as we will see in Section~\ref{sec:bifurcation}.

As for the phase interaction function in~\eqref{eq:main_system}, it only depends on the first harmonic of the state of the oscillator at~$x$. Hence, the system can be reduced to the Ott--Antonsen manifold~\cite{Ott2008,Bick2018c} as for globally coupled networks~\cite{Skardal2019a,Skardal2019b}.

\section{Bifurcations of Twisted States}
\label{sec:bifurcation}
\noindent
In this section, we study the bifurcation that occurs when a $q$-twisted state gains or loses its stability under variation of parameters. 
As noted above, the system~\eqref{eq:main_system} has a $\T$-symmetry that maps a solution $\theta_x$ to a solution $\theta_x+\beta$ for a given constant $\beta\in\T$. 
Thus, the linearization of the right-hand side of~\eqref{eq:main_system} always has a zero eigenvalue, which makes a rigorous bifurcation analysis tedious. 
In order to avoid this zero eigenvalue, we change to the system of phase differences and define $\Psi_x(t) := \theta_x(t)-\theta_0(t)$. Then, the function $\Psi_x(t)$ satisfies
\begin{subequations}\label{eq:main_pd}
\begin{align}
\begin{split}\label{eq:main_pd_pairwise}
    \dot \Psi_x &= K_2 \int_\S G(x-y) \sin(\Psi_y-\Psi_x+\alpha_2)\ \d y\\
    &\qquad - K_2\int_\S G(-y) \sin(\Psi_y + \alpha_2)\ \d y
\end{split}\\
\begin{split}\label{eq:main_pd_triplet}
    &\quad\ + K_3\int_\S \int_\S G(x-y)G(x-z) \sin(2\Psi_y - \Psi_z-\Psi_x + \alpha_3)\ \d z \d y\\
    &\qquad - K_3 \int_\S\int_\S G(-y)G(-z) \sin(2\Psi_y - \Psi_z + \alpha_3) \ \d z \d y,
\end{split}
\end{align}
\end{subequations}
and 
\begin{align}\label{eq:pd_boundary_cond}
    \Psi_0(t) = 0
\end{align}
 for all times $t$. The process of transitioning from \eqref{eq:main_system} to the system of phase differences \eqref{eq:main_pd} reduces the $\T$-symmetry. In particular, the $\T$-symmetry is present in the system \eqref{eq:main_system} but not in \eqref{eq:main_pd}, since every function in the system of phase differences has to satisfy $\Psi_0(t) = 0$ and thus cannot be shifted by a constant. After this symmetry reduction, the $q$-twisted states \eqref{eq:SplayContinuum} are represented by the function $\Psi^q_x = \Psi^q(x)$ with
 \begin{align}\label{eq:twisted_pd}
     \Psi^q(x) := 2\pi q x,
 \end{align}
 which does not depend on the time $t$, satisfies $\Psi^q(0)=0$, and cannot be perturbed by a constant function. Consequently, when linearizing the right-hand side of \eqref{eq:main_pd}, which we denote by $\mathcal G$, around a $q$-twisted state \eqref{eq:twisted_pd} there is no trivial zero eigenvalue. To make this precise, we define a function $\mathcal G^q \colon X \times \mathcal P\to X$ as
\begin{align*}
    \mathcal G^q(v,p) = \mathcal G(\Psi^q + v, p),
\end{align*}
where $q\in\Z$ is the winding number, $p = (K_2, K_3, A, B, \alpha_2, \alpha_3)$ summarizes all parameters and $ X = H_0^1(\S,\R)$ is the space of once weakly differentiable functions on $\S$ with zero boundary conditions. Since $H^1(\S,\R)\subset C(\S)$, these boundary conditions can be imposed in the classical sense. The function~$\mathcal G^q$ consequently gives the local behavior of the right-hand side of \eqref{eq:main_pd} around the $q$-twisted state, $v\in X$ can be seen as a perturbation of the twisted state, and the condition $v(0)=0$ ensures~\eqref{eq:pd_boundary_cond}. Conducting a bifurcation analysis of~\eqref{eq:main_pd} instead of~\eqref{eq:main_system} simplifies the setting.




\subsection{Linearization}
\label{sec:Lin}

The bifurcation analysis of twisted states is based on a stability analysis, which can be conducted using the eigenvalues of the linearization of the right-hand side around the twisted state. More precisely, we linearize $\mathcal G^q(v,p)$ around $v=0$, consider this linearization as an operator $D_v \mathcal G^q(0,p)$ from $X$ to itself and determine the eigenvalues of this operator.

To obtain $D_v \mathcal G^q(0,p)$, we take a function $\eta_\cdot \in X$ and $h\in\R$, and calculate
\begin{align*}
    D_v \mathcal G^q(0,p)[\eta] = \lim_{h\to 0}\frac{1}{h} \mathcal G^q(h\eta, p).
\end{align*}
Using the definition of $\mathcal G^q$ we obtain
\begin{align*}
    D_v\mathcal G^q(0,p)[\eta] &= K_2\int_\S G(x-y) (\eta_y - \eta_x)\cos(\Psi^q_y-\Psi^q_x + \alpha_2)\ \d y\\
    &\qquad - K_2 \int_\S G(-y) \eta_y \cos(\Psi^q_y+\alpha_2)\ \d y\\
    &\quad\ + K_3 \int_\S \int_\S G(x-y)G(x-z)(2\eta_y-\eta_z+\eta_x)\\
    &\qquad\qquad\qquad\cdot\cos(2\Psi^q_y - \Psi^q_z - \Psi^q_x+\alpha_3) \ \d z \d y\\
    &\qquad - K_3 \int_\S \int_\S G(-y)G(-z)(2\eta_y-\eta_z)\\
    &\qquad\qquad\qquad\quad\cdot\cos(2\Psi^q_y+\Psi^q_z+\alpha_3) \ \d z \d y.
\end{align*}
Even though this computation was very formal, it can be rigorously shown that $D_v\mathcal G^q(0,p)$, as calculated here, is indeed the Fr\'echet derivative of $\mathcal G^q(v,p)$ at $v=0$, see~\cite{Bick2023a}.

Henceforth we focus on the twisted state with the smallest winding number, namely $q=1$; the case $q=-1$ is analogous.
We note that the functions $u_k(x) = \sin(2\pi k x)$ and $w_k(x) = 1-\cos(2\pi k x)$ for $k\ge 1$ form a Schauder basis of~$X$, see~\cite[Lemma B.1]{Bick2023a}. 
Consequently, we evaluate $D_v\mathcal G^1(0,p_0)$ on these basis functions. For $k=1$, this yields
\begin{align*}
    D_v \mathcal G^1(0,p)[u_1] &= \frac {K_2} 2 \Big(\cos(\alpha_2)-A\cos(\alpha_2) - B \sin(\alpha_2)\Big)u_1 + \frac {K_2} 2 \sin(\alpha_2)w_1\\
    &\quad + \frac{K_3}4(A^2+B^2)\cos(\alpha_3) u_1 + \frac{K_3}4(A^2+B^2)\sin(\alpha_3)w_1,\\
    D_v \mathcal G^1(0,p)[w_1] &= -\frac{K_2}2\sin(\alpha_2)u_1 + \frac{K_2}2\Big(\cos(\alpha_2)-A\cos(\alpha_2)-B\sin(\alpha_2)\Big)w_1\\
    &\quad - \frac{K_3}{4} (A^2+B^2)\sin(\alpha_3) u_1 + \frac{K_3}4 (A^2+B^2)\cos(\alpha_3)w_1,
\end{align*}
which implies a complex conjugated pair of eigenvalues
\begin{subequations}\label{eq:OneTwLinStab}
\begin{align}
\begin{split}
    \lambda_1^\pm &= \frac{K_2}{2}\left( \cos(\alpha_2) - A\cos(\alpha_2)-B\sin(\alpha_2) \pm i \sin(\alpha_2)\right)\\
	&\quad + \frac{K_3}{4}\left( (A^2+B^2)\cos(\alpha_3) \pm  i (A^2+B^2)\sin(\alpha_3) \right),
	\end{split}
\intertext{where $i = \sqrt{-1}$. Similarly, one can insert the basis functions $u_k, w_k$ for $k\ge 2$ into $D_v \mathcal G^1(0,p_0)$, which results in further eigenvalues}
\begin{split}
	\lambda_2^\pm &= \frac{K_2}{4}\left(-A\cos(\alpha_2) - 3B\sin(\alpha_2) \pm i(B\cos(\alpha_2) + A\sin(\alpha_2))\right)\\
	&\quad\ +\frac{K_3}4\left(2A\cos(\alpha_3)-2B\sin(\alpha_3)\pm i (2B\cos(\alpha_3)+2A\sin(\alpha_3))\right) 
	\end{split}\\
	\begin{split}
	\lambda_3^\pm &= -\frac {K_2}2 (A\cos(\alpha_2) + B\sin(\alpha_2))\\
	&\quad\ + \frac {K_3}4 ((A^2-B^2)\cos(\alpha_3) - 2AB\sin(\alpha_3))\\
	&\quad\ \pm \frac {K_3}4 i ((A^2-B^2)\sin(\alpha_3) + 2AB \cos(\alpha_3))
	\end{split}\\
	\lambda_k^\pm &=  -\frac {K_2}2 (A\cos(\alpha_2) + B\sin(\alpha_2))
\end{align}
\end{subequations}
for $k\ge 4$. Moreover, corresponding eigenfunctions are given by $v_k^\pm = w_k \pm i u_k$ for $k\in\N$.



\begin{figure}
    \centering
    \begin{overpic}[width = 0.9\textwidth]{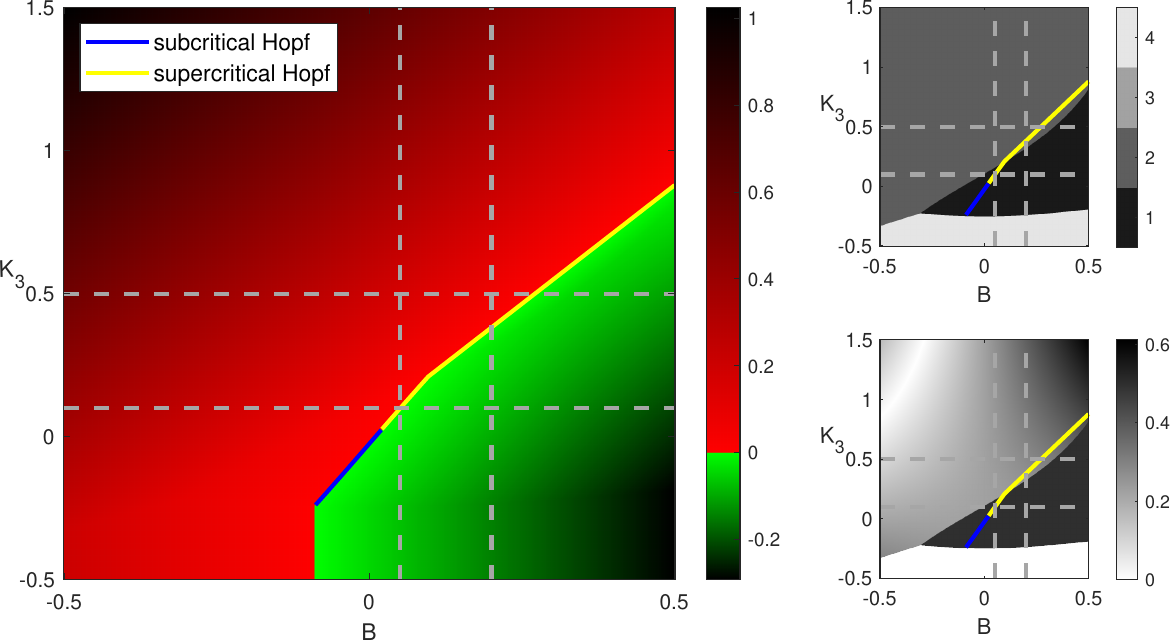}
        \put(0,0){{\small\textbf{(a)}}}
        \put(70,29){\small\textbf{(b)}}
        \put(70,0){\small\textbf{(c)}}
    \end{overpic}
    \caption{Eigenvalues of the linearization of a $1$-twisted state. Panel~(a) shows $\max_k \Real(\lambda_k^+)$. In the green regions, the $1$-twisted state is stable, whereas it is unstable in the red regions. At the boundary a bifurcation occurs. Panel~(b) depicts which eigenvalue is the critical one, i.e., for which~$\ell$ we have $\Real(\lambda_\ell^+) = \max_k \Real(\lambda_k^+)$. Finally, Panel~(c) shows $|\Imag(\lambda_\ell^+)|$, i.e., the modulus of the imaginary part of the critical eigenvalue. Parameter values: $K_2 = 1, A=0.9, \alpha_2 = \frac{\pi}{2} - 0.1, \alpha_3 = 0$.
    Dashed lines in all panels indicate the parameter values for the numerical continuation shown in Figs.~\ref{Fig:ScanB_K0_1}--\ref{Fig:ScanK_B0_2} below.}
    \label{fig:Eigenvalues}
\end{figure}

Whenever $\Imag(\lambda_\ell^+)\neq 0$ but $\Real(\lambda_\ell^+)=0$, for some $\ell$, we expect a Hopf bifurcation to occur under variation of parameters. The transverse stability of this bifurcation is then determined by the other eigenvalues. If there exists $k\in\N$ with $k\neq \ell$ such that $\Real(\lambda_k^+)>0$, all equilibria and periodic orbits that emanate from the Hopf bifurcation are not transversely stable. If on the other hand $\Real(\lambda_k^+)<0$ for all $k\neq \ell$, the equilibria and periodic orbits are at least transversely stable. As we are particularly interested in transversely stable bifurcations, we assume from now on that this is the case.
As one can see in Figure \ref{fig:Eigenvalues}, there are parameter regions where the critical eigenvalue is attained for $\ell = 1$ and other regions where $\ell=2$ is the dominant eigenvalue. Here, a Hopf bifurcation can occur. However, there also exist parameter regions, where $\ell=4$ is the critical eigenvalue. Since $\Imag(\lambda_4^\pm) = 0$, there is no generic Hopf bifurcation for these parameter values.

\subsection{Bifurcations}\label{sec:bif}
Having investigated the linear stability of the $1$-twisted state, we now analyze the periodic solutions that emanate from the $1$-twisted state in a Hopf bifurcation.

We use the notation $C_T(\R, X)$ for all continuous $T$-periodic functions with values in~$X$. Moreover, we assume for simplicity that~$B$ is the main bifurcation parameter, i.e., we fix all other parameters and vary only~$B$ to initiate the bifurcation, which then occurs at some value~$B_0$. We can then write the eigenvalues of linearization around the $1$-twisted state as functions of~$B$.
At the critical value~$B_0$ we have $\Real(\lambda_\ell^\pm(B_0)) = 0$ and we additionally assume that $\Real(\lambda_k^\pm(B_0))<0$ for all $k\neq \ell$ and $\Imag(\lambda_\ell^\pm(B_0))\neq 0$. In particular, we also denote $\lambda_\ell(B)$ for the critical eigenvalue that has positive imaginary part, i.e., $\lambda_\ell(B) := \lambda_\ell^+(B)$ if $\Imag(\lambda_\ell^+(B_0))>0$ and $\lambda_\ell(B) := \lambda_\ell^-(B)$ else. Similarly, we denote $v_\ell$ for the eigenfunction that corresponds to this critical eigenvalue.

A generic theorem, that relies on some technical assumptions \cite[Theorem I.8.2]{Kielhofer2012}, then guarantees that periodic solutions bifurcate from the $1$-twisted state when~$B$ passes through~$B_0$. In particular, there is a continuously differentiable function $\kappa\colon U\to \R$, where $U$ is a small neighborhood of $0\in \R$, such that $\Imag(\lambda_\ell(B_0))=\kappa(0)>0$. Moreover, there is a continuous curve $U\ni r\mapsto (B(r), v(r))\in (B_0-\delta, B_0+\delta)\times C_{2\pi/\kappa(r)}(\R, X)$, such that for every $r\in U$ the function $\Psi^1 + v(r)$ is a solution of \eqref{eq:main_pd} when the parameter is set to $B(r)$. This curve satisfies $B(0) = B_0$ and $v(0) \equiv 0\in C_{2\pi/\kappa(0)}(\R,X)$.
Furthermore, every other periodic solution of \eqref{eq:main_system} in a neighborhood of the $1$-twisted state for parameter values $B(r)$ can be obtained as a phase shift from $\Psi^1 + v(r)$, i.e., it is given by $\Psi^1 + S_\tau v(r)$, where $(S_\tau v(r))(t) := v(r)(t+\tau)$, see~\cite{Kielhofer2012}.


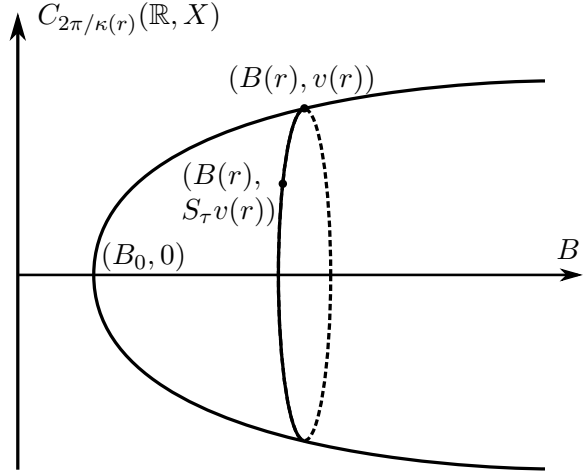
\begin{figure}
    \centering
    \begin{overpic}[width = .6\textwidth]{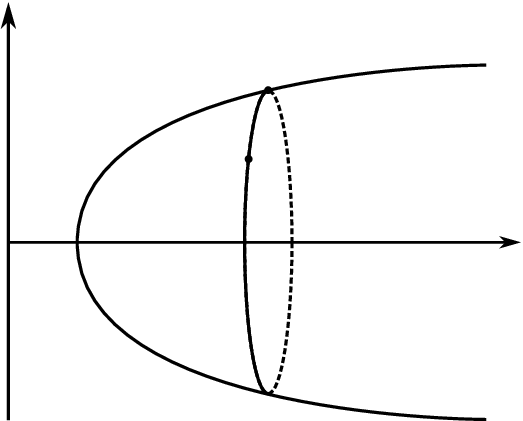}
        \put(95,37){$B$}
        \put(5,78){$C_{2\pi/\kappa(r)}(\R, X)$}
        \put(16,36){$(B_0,0)$}
        \put(38,67){$(B(r), v(r))$}
        \put(30,50){$(B(r),$}
        \put(30,44){$S_\tau v(r))$}
    \end{overpic}
    \caption{Illustration of the Hopf bifurcation. Note that the vertical axis represents the space $C_{2\pi/\kappa(r)}$, which is different from classical bifurcation diagrams. Thus, there is also no flow around the paraboloid, but each point on the paraboloid represents a solution. Each point on the black parabola can be reached via the curve $r\mapsto (B(r), v(r))$. Every other periodic solution can be obtained by a phase shift $(B(r), S_\tau v(r))$ for a suitable $\tau\in\R$.}
    \label{fig:BifurcationSketch}
\end{figure}

Additionally to the first derivatives of~$\mathcal G^1$, second and third derivatives are required to approximate the curve $(v(r), B(r))$, see~\cite[Theorem I.9.1]{Kielhofer2012}. In particular, the curve of $2\pi/\kappa(r)$-periodic solutions can be approximated as
\begin{subequations}\label{eq:curve_approx}
\begin{align}
    \label{eq:amplitude_approx}
    \frac{\d}{\d r} v(r) \Big\vert_{r=0}(t) &= 2\Real(v_\ell e^{ i \kappa(0) t}),\\
    \label{eq:period_approx}
    \frac{\d}{\d r} \kappa(r)\vert_{r=0} &= 0,\\
    \label{eq:param_approx}
    \frac{\d^2}{\d r^2}B(r)\Big\vert_{r=0} &= \frac{\Real(\zeta)}{\Real(\frac{\d}{\d B}\lambda_\ell(B)\vert_{B=B_0})},
\end{align}
\end{subequations}
where $\zeta\in \mathbb C$ can be computed from the first, second and third derivative of the right-hand side~$\mathcal G$ at the $1$-twisted state; see Appendix~\ref{sec:bifurcation_formulas}. 
Further, $\Real(\frac{\d}{\d B}\lambda_\ell(B)\vert_{B=B_0})$ is the speed with which the real part of the critical eigenvalue passes through zero. 
Basically,~\eqref{eq:amplitude_approx} helps to approximate the amplitude and profile of bifurcating periodic solutions, \eqref{eq:period_approx}~determines their period, and \eqref{eq:param_approx}~connects the parameter~$B$ with~$r$ and thereby determines for which parameter value of~$B$ these periodic solutions occur.
In particular, if~\eqref{eq:param_approx} is positive, these periodic solutions are only existent when $B>B_0$. 
Conversely, if~\eqref{eq:param_approx} is negative, they exist for $B<B_0$. Using the principle of exchange of stability of the equilibrium with the periodic orbits in a Hopf bifurcation, one can even determine the stability of the periodic orbits. Specifically, if $\Real(\zeta)>0$ one can consider two cases depending on the sign of $\Real(\frac{\d}{\d B}\lambda_\ell(B)\vert_{B=B_0})$. 
If it is positive, then the $1$-twisted state is unstable for $B>B_0$. 
Moreover, \eqref{eq:param_approx}~is positive and thus the periodic solutions only exist when $B>B_0$. 
Since the twisted states are unstable in that parameter regime, the periodic solutions have to be stable. 
If on the other hand $\Real(\frac{\d}{\d B}\lambda_\ell(B)\vert_{B=B_0})$ is negative, the twisted state is unstable for $B<B_0$. 
Since~\eqref{eq:param_approx} is then also negative, the periodic solutions also exist in the same parameter region and are stable. 
In both cases the bifurcation is supercritical. Repeating the argument for $\Real(\zeta)<0$ shows that the bifurcation is then subcritical, meaning that the periodic solutions that emanate from the bifurcation are then unstable. 
In conclusion, the sign of $\Real(\zeta)$ can be used to distinguish between sub- and supercritical Hopf bifurcations.

\begin{rem}\label{rem:Period}
%
Note that the index~$\ell$ of the critical eigenvalue determines the periodicity of the eigenfunctions~$v^\pm_\ell$, which in turn relates to the periodicity of the (linear approximation) of the bifurcating solutions.
More specifically, let us consider the Hopf bifurcation
corresponding to the critical eigenvalue $\lambda_\ell^+$ with the eigenfunction $v_\ell^+ = w_\ell + i u_\ell = 1 - e^{-2 \pi \ell i x}$. 
In this case, formula~\eqref{eq:amplitude_approx} gives
\[v(r) \approx 2 r \Real( v_\ell^+ e^{i \kappa(0) t} )
= 2 r \Real( ( 1 - e^{-2 \pi \ell i x} ) e^{i \kappa(0) t} ).\]
Hence, for $B\approx B_0$ equation~\eqref{eq:main_pd} has a solution curve with an $r$-parametric asymptotic representation
\[\Psi_x(t) \approx 2 \pi x + 2 r \Real( e^{i \kappa(0) t} ) - 2 r \Real( e^{-2 \pi \ell i ( x - \kappa(0) t / (2 \pi \ell) )} ).\]
Recalling that $\Psi_x(t) = \theta_x(t) - \theta_0(t)$,
where~$\theta_x(t)$ is a solution of equation~\eqref{eq:main_system}, we conclude that
\[\theta_x(t) \approx 2 \pi x -
2 r \cos( 2 \pi \ell ( x - \kappa(0) t / (2 \pi \ell) ) ).\]
In other words, in the new periodic solution of equation~\eqref{eq:main_system}, we expect a small-amplitude perturbation on top of the $1$-twisted state with the phase profile~$\cos( 2 \pi \ell x)$ drifting with the speed $s = \kappa(0) / (2 \pi \ell) = \Imag(\lambda_\ell^+) / (2 \pi \ell)$.

\end{rem}

\section{Continuation on the Ott--Antonsen Manifold}
\label{sec:OttAnt}
\noindent
It is easy to verify that the dynamics described by Eq.~\eqref{eq:main_system} is equivalent to the dynamics of the Ott--Antonsen equation
\begin{equation}
\begin{split}
\pf{}{t} z(x,t) &= \fr{K_2}{2}\left(e^{i \alpha_2} \mathcal{G} z - e^{-i \alpha_2} z^2 \mathcal{G} \overline{z}\right)\\&
\qquad + \fr{K_3}{2}\left(e^{i \alpha_3} \mathcal{G} z^2\:\mathcal{G} \overline{z} - e^{-i \alpha_3} z^2 \mathcal{G} \overline{z}^2\:\mathcal{G} z\right),
\end{split}
\label{Eq:OA}
\end{equation}
with a convolution-type integral operator
\[
\mathcal{G} z = \int_{y\in\mathbb{S}} G(x - y) z(y,t) \d y,
\]
after insertion of the ansatz
\[
z(x,t) = e^{i \Theta(x,t)}.
\]
Eq.~\eqref{Eq:OA} is useful for two reasons.
First, it can be used to perform a linear stability analysis of twisted states in an alternative way; see Appendix~\ref{Appendix:OA}.
Second, this equation can be used to compute the solution branches emanating from the Hopf bifurcations found above.
In this section, we outline how to compute these solution branches; the approach is based on recent results presented in~\cite{Omelchenko2022} and adapted to our setting.

\subsection{Self-consistency equation for traveling nonuniformly twisted waves}

All periodic solutions emerging at the Hopf bifurcations of $1$-twisted states
have form of quasiperiodically evolving synchronization patterns
\[
\Theta(x,t) = \Omega t + \Theta_0(x - s t)
\]
with some~$\Theta_0(x)$ and $s,\Omega\in\mathbb{R}$.
For Eq.~\eqref{Eq:OA} they correspond to traveling wave solutions of the form
\begin{equation}
z(x,t) = a(x - s t) e^{i \Omega t},
\label{Ansatz:Wave}
\end{equation}
where $|a(x)| = 1$ for all $x\in\mathbb{S}$.
Such solutions can be efficiently computed, using the self-consistency equation derived below.

By inserting ansatz~\eqref{Ansatz:Wave} into Eq.~\eqref{Eq:OA} we find that the profile function~$a(x)$ in~\eqref{Ansatz:Wave} is a periodic solution of the integro-differential equation
\begin{equation}
    \begin{split}
-sa'(x) &= - i \Omega a + \fr{K_2}{2} e^{i \alpha_2} \mathcal{G} a - \fr{K_2}{2} e^{-i \alpha_2} a^2 \mathcal{G} \overline{a}\\ 
&\qquad+ \fr{K_3}{2} e^{i \alpha_3} \mathcal{G} a^2\:\mathcal{G} \overline{a} - \fr{K_3}{2} e^{-i \alpha_3} a^2 \mathcal{G} \overline{a}^2\:\mathcal{G} a.
    \end{split}
    \label{Eq:a}
\end{equation}
If $s\ne 0$, Eq.~\eqref{Eq:a} is equivalent to the complex Riccati equation
\begin{equation}
a'(x) = w(x) + i \frac{\Omega}{s} a - \overline{w}(x) a^2,
\label{Eq:Riccati}
\end{equation}
with
\begin{equation}
w(x) = - \fr{K_2}{2 s} e^{i \alpha_2} \mathcal{G} a - \fr{K_3}{2 s} e^{i \alpha_3} \mathcal{G} a^2\:\mathcal{G} \overline{a}
\label{Def:w}
\end{equation}
contains the network coupling terms.
Rather than finding the unknown function~$a(x)$, our strategy is to determine the corresponding mean field~$w(x)$ using a self-consistency equation.

First, we recall some facts about Eq.~\eqref{Eq:Riccati}, which were previously proved in~\cite[Section~2]{Omelchenko2022} by one of the authors of this paper.
The most important fact is that for an arbitrary periodic function~$w(x)$ and an arbitrary real coefficient~$\Omega / s$, the complex Riccati equation~\eqref{Eq:Riccati} usually has two periodic solutions.
The initial conditions of these solutions
are determined by the fixed points
of the Poincar{\'e} map of Eq.~\eqref{Eq:Riccati}.
Moreover, due to the special structure of Eq.~\eqref{Eq:Riccati}, its Poincar{\'e} map
has the form of a M{\"o}bius transformation
\[
\mathcal{M}(z) = \frac{e^{i\theta}(z + b)}{\overline{b}z + 1}
\]
with some $b\in\{z\in\mathbb{C}\::\: |z|<1 \}$ and $\theta\in 2\pi\mathbb{S}$ determined by the choice of~$w(x)$ and~$\Omega / s$.
If $|b| > |\sin(\theta/2)|$, then the both fixed points of~$\mathcal{M}(z)$ lie on the unit circle and therefore Eq.~\eqref{Eq:Riccati} has two periodic solutions satisfying $|a(x)| = 1$.
One of these solutions is asymptotically stable, while the other is asymptotically unstable.
In contrast, if $|b| < |\sin(\theta/2)|$, then the fixed points of~$\mathcal{M}(z)$ lie inside or outside of the unit circle.
More specifically, in this case, Eq.~\eqref{Eq:Riccati} has one solution that satisfies the inequality $|a(x)| < 1$ and another solution that satisfies the inequality $|a(x)| > 1$.

Suppose that the Ott--Antonsen equation~\eqref{Eq:OA} has a stable traveling wave solution of the form~\eqref{Ansatz:Wave} with $|a(x)| = 1$.
For this solution, we can calculate its mean field~$w(x)$ through~\eqref{Def:w}.
Then, considering~$w(x)$ and~$\Omega / s$ as given, we can try to solve the periodic boundary value problem for Eq.~\eqref{Eq:Riccati} with the constraint $|a(x)| = 1$.
Above we have shown that such a problem can either have two solutions (one stable and one unstable) or none.
Moreover, in the first case, depending on the sign of the speed~$s$, only one of the two solutions can be relevant for a stable traveling wave~\eqref{Ansatz:Wave}.
Indeed, recalling that in~\eqref{Ansatz:Wave} we use the moving coordinate frame $\xi = x - s t$, we easily conclude that if $s > 0$ ($s < 0$) then a stable traveling wave~(\ref{Ansatz:Wave}) corresponds to an unstable (stable) solution of Eq.~\eqref{Eq:Riccati}.
Altogether these facts allow us to say that using the periodic boundary value problem for Eq.~\eqref{Eq:Riccati} we can always reconstruct
the profile~$a(x)$ of a stable traveling wave~\eqref{Ansatz:Wave} if the corresponding
mean field~$w(x)$ and the ratio~$\Omega / s$ are known.
Denoting the resulting solution operator by $\mathcal{U}$ we can write
\begin{equation}
a(x) = \mathcal{U}( w(x), \Omega/s ).
\label{Operator:U}
\end{equation}
Obviously, the last expression agrees with the definition of~$w(x)$ if and only if
\begin{equation}
    \begin{split}
w(x) &= - \fr{K_2}{2 s} e^{i \alpha_2} \mathcal{G} \mathcal{U}( w(x), \Omega/s ) \\ 
&\qquad - \fr{K_3}{2 s} e^{i \alpha_3} \mathcal{G} [\mathcal{U}( w(x), \Omega/s )]^2\:\mathcal{G} \overline{\mathcal{U}( w(x), \Omega/s )},
    \end{split}
    \label{Eq:SC}
\end{equation}
which is an integral \emph{self-consistency equation} for~$w(x)$.

To complete the definition of Eq.~(\ref{Eq:SC}), we show a simple and fast way to calculate the operator~$\mathcal{U}$.
The justification of this method, consisting of five steps, is given in~\cite[Section~2]{Omelchenko2022}:

(i)~Given~$w(x)$ and~$\Omega / s$, solve Eq.~\eqref{Eq:Riccati} starting from the initial condition $a(0) = 1$, and denote $\zeta_1 = a(1)$.

(ii)~Similarly, solve Eq.~\eqref{Eq:Riccati} in the backward time starting from the initial condition $a(0) = 0$, and denote $\zeta_0 = a(-1)$.

(iii)~Calculate the coefficients~$b$ and~$e^{i\theta}$ of the M{\"o}bius transformation~$\mathcal{M}(z)$ representing the Poincar{\'e} map of Eq.~\eqref{Eq:Riccati},
\[
b = - \zeta_0\qquad\text{and}\qquad
e^{i\theta} = \frac{\overline{\zeta}_0 - 1}{\zeta_0 - 1} \zeta_1.
\]
(Note that by construction $|\zeta_1| = 1$.)

(iv) Check that $|b| > |\sin(\theta/2)|$ (otherwise the operator $\mathcal{U}$ is not well-defined) and then calculate the initial value~$a_*$ of the periodic solution of interest, namely
\begin{align*}
a_* &= \frac{i \sin(\theta/2) - \sqrt{ |b|^2 - \sin^2(\theta/2) }}{|b|^2} b e^{i \theta/2}\quad\text{if $s > 0$},\\
\intertext{or}
a_* &= \frac{i \sin(\theta/2) + \sqrt{ |b|^2 - \sin^2(\theta/2) }}{|b|^2} b e^{i \theta/2}\quad\text{if $s < 0$}.
\end{align*}

(v)~Solve Eq.~\eqref{Eq:Riccati}, starting from the initial condition $a(0) = a_*$. This yields the periodic solution of Eq.~\eqref{Eq:Riccati}.

\begin{rem}
If in step~(iv) of the above algorithm we choose the formula of~$a_*$ with~`$-$' for $s < 0$ and the formula with~`$+$' for $s > 0$, we obtain another solution operator for the periodic boundary value problem associated with Eq.~\eqref{Eq:Riccati}.
But this operator, by construction, gives profile functions~$a(x)$ of traveling waves~\eqref{Ansatz:Wave}, which are unstable with respect to Eq.~\eqref{Eq:OA}.
\end{rem}

\subsection{Algebraic self-consistency equation}

The self-consistency equation~\eqref{Eq:SC} is a nonlinear integral equation that can be difficult to solve.
However, in the case of coupling kernel~\eqref{Formula:G}, it can be reduced to a finite-dimensional nonlinear system. To see this, let us denote
\begin{align*}
\psi_1(x) &= 1, & \psi_2(x) &= \cos(2\pi x), & \psi_3(x) &= \sin(2\pi x),\\
&&\psi_4(x) &= \cos(4\pi x), & \psi_5(x) &= \sin(4\pi x).
\end{align*}
Then using trigonometric identities we can write
\[
\mathcal{G} u = \langle u, \psi_1 \rangle \psi_1
+ ( A  \langle u, \psi_2 \rangle - B \langle u, \psi_3 \rangle ) \psi_2
+ ( A \langle u, \psi_3 \rangle + B \langle u, \psi_2 \rangle ) \psi_3,
\]
where
\[
\langle u, v \rangle = \int_{x\in\mathbb{S}} u(x) \overline{v}(x) \d x
\]
is the usual~$L^2$ inner product. 
Moreover, for complex numbers $a_0$, $a_1$, $a_2$, $b_0$, $b_1$, $b_2$ we have
\begin{align*}
&( a_1 \psi_1 + a_2 \psi_2 + a_3 \psi_3 ) ( b_1 \psi_1 + b_2 \psi_2 + b_3 \psi_3 ) \\
&\qquad
= \left[ a_1 b_1 + \fr{1}{2} ( a_2 b_2 + a_3 b_3 ) \right] \psi_1 + ( a_2 b_1 + a_1 b_2 ) \psi_2 \\
&\qquad\qquad
+ ( a_3 b_1 + a_1 b_3 ) \psi_3 + \fr{1}{2} ( a_2 b_2 - a_3 b_3 ) \psi_4 + \fr{1}{2} ( a_3 b_2 + a_2 b_3 ) \psi_5.
\end{align*}
This relation together with Eq.~\eqref{Eq:SC} implies
\begin{equation}
w(x) = \sum\limits_{k=1}^5 \hat{w}_k \psi_k(x)
\label{Ansatz:w}
\end{equation}
with some complex coefficients~$\hat{w}_k$.
Inserting ansatz~\eqref{Ansatz:w} into Eq.~\eqref{Eq:SC} and equating the terms proportional to different~$\psi_k$ separately,
we reformulate Eq.~\eqref{Eq:SC} as a system of five complex equations
\begin{equation}
\begin{split}
    - 2 s \hat{w}_1 &=  K_2 e^{i \alpha_2} a_1 +  K_3 e^{i \alpha_3} \left( \overline{a}_1 b_1 + \fr{1}{2} ( \overline{a}_2 b_2 + \overline{a}_3 b_3 ) \right),\\
- 2 s \hat{w}_2 &=  K_2 e^{i \alpha_2} a_2 +  K_3 e^{i \alpha_3} \left( \overline{a}_2 b_1 + \overline{a}_1 b_2 \right),\\
- 2 s \hat{w}_3 &=  K_2 e^{i \alpha_2} a_3 +  K_3 e^{i \alpha_3} \left( \overline{a}_3 b_1 + \overline{a}_1 b_3 \right),\\
- 2 s \hat{w}_4 &=  \fr{K_3}{2} e^{i \alpha_3} \left( \overline{a}_2 b_2 - \overline{a}_3 b_3 \right),\\
- 2 s \hat{w}_5 &=  \fr{K_3}{2} e^{i \alpha_3} \left( \overline{a}_3 b_2 + \overline{a}_2 b_3 \right),
\end{split}
    \label{System:w}
\end{equation}
with
\begin{align*}
a_1 &= \Psi_1, & a_2 &= A \Psi_2 - B \Psi_3, & a_3 &= A \Psi_3 + B \Psi_2,\\
b_1 &= \Phi_1, & b_2 &= A \Phi_2 - B \Phi_3,& b_3 &= A \Phi_3 + B \Phi_2,
\end{align*}
and
\begin{align*}
\Psi_k &= \left\langle \mathcal{U}\left( \sum\limits_{k=1}^5 \hat{w}_j \psi_j, \frac{\Omega}{s} \right), \psi_k \right\rangle,\\
\Phi_k &= \left\langle \left[ \mathcal{U}\left( \sum\limits_{k=1}^5 \hat{w}_j \psi_j, \frac{\Omega}{s}\right) \right]^2, \psi_k \right\rangle .
\end{align*}
System~\eqref{System:w} needs to be solved with respect to~$s$, $\Omega$ and $\hat{w}_k$, $k=1,\dots,5$.
But it is obviously underdetermined for this.
The problem can be resolved by recalling that the original Eq.~\eqref{Eq:OA} has two continuous symmetries. Therefore, we may add two pinning conditions to~\eqref{System:w}.
For the sake of convenience, we choose these pinning conditions in the form
\begin{equation}
\Imag \hat{w}_2 = \Real \hat{w}_3 = 0.
\label{Pinning:w}
\end{equation}
Now, if we find a solution of the extended system~\eqref{System:w}, \eqref{Pinning:w}, we can use formula~\eqref{Ansatz:w} to calculate~$w(x)$ and then~\eqref{Operator:U} to calculate the corresponding~$a(x)$.
Altogether, two scalars~$s$ and~$\Omega$ and the profile function~$a(x)$ allow us to determine the traveling wave solution~\eqref{Ansatz:Wave} of Eq.~\eqref{Eq:OA}.
Finally, if we want to show a typical snapshot of the corresponding solution of Eq.~\eqref{eq:main_system}, we can use $\Theta(x) = \mathrm{arg}\:a(x)$.

\section{Continuation of Periodic Traveling Solutions}
\label{sec:Continuation}
\begin{figure}
\begin{center}
\includegraphics[width=0.64\textwidth]{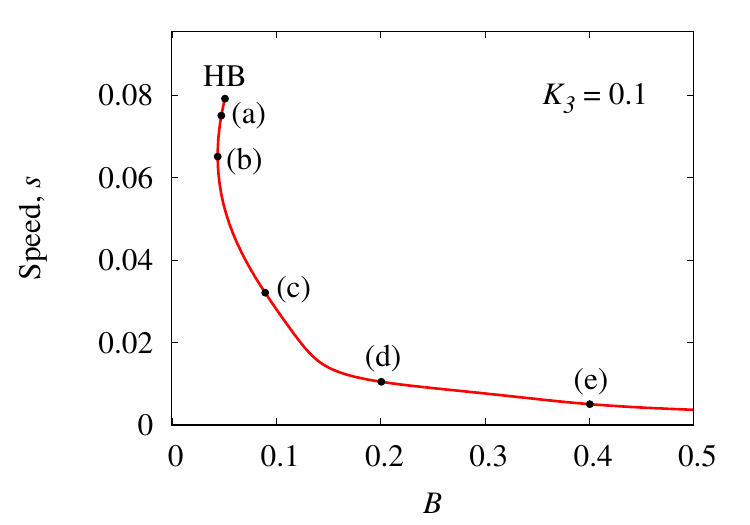}
\begin{minipage}[b]{0.32\textwidth}
\includegraphics[width=1.0\textwidth]{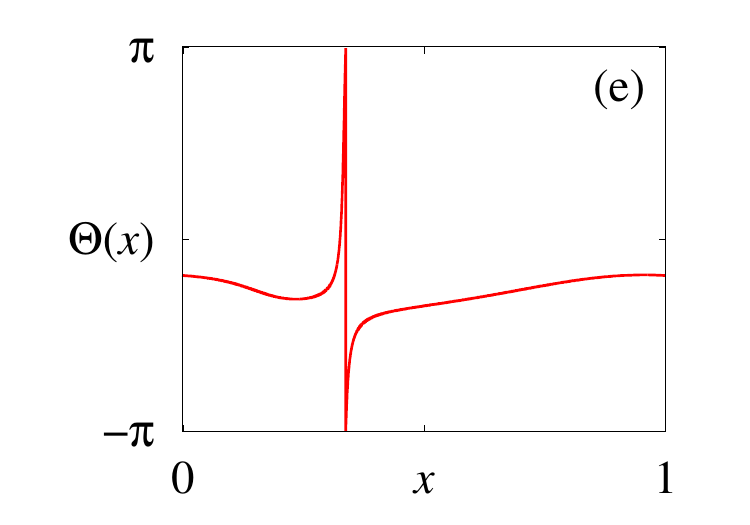}\\[2mm]
\includegraphics[width=1.0\textwidth]{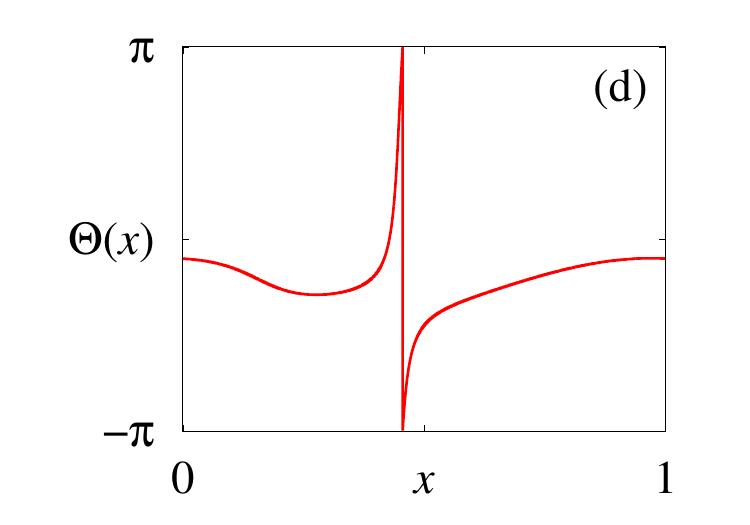}
\end{minipage}
\\[2mm]
\includegraphics[width=0.32\textwidth]{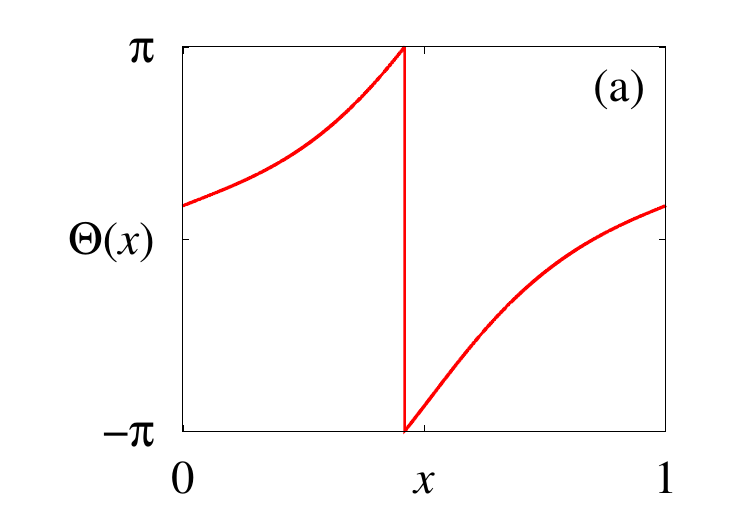}
\includegraphics[width=0.32\textwidth]{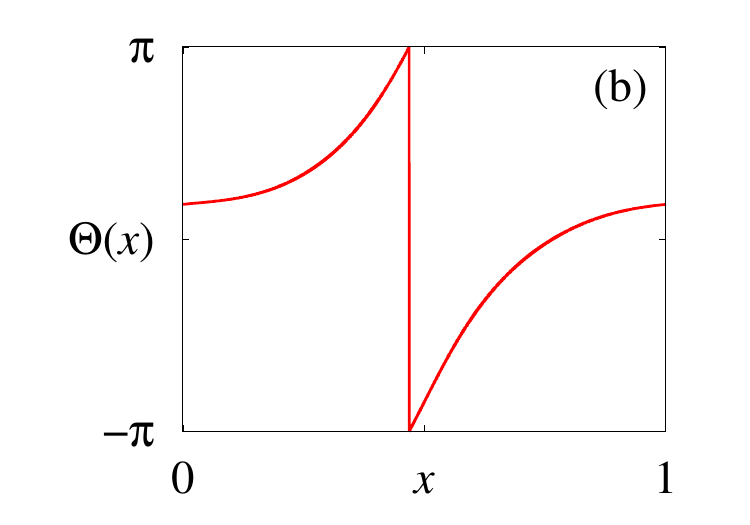}
\includegraphics[width=0.32\textwidth]{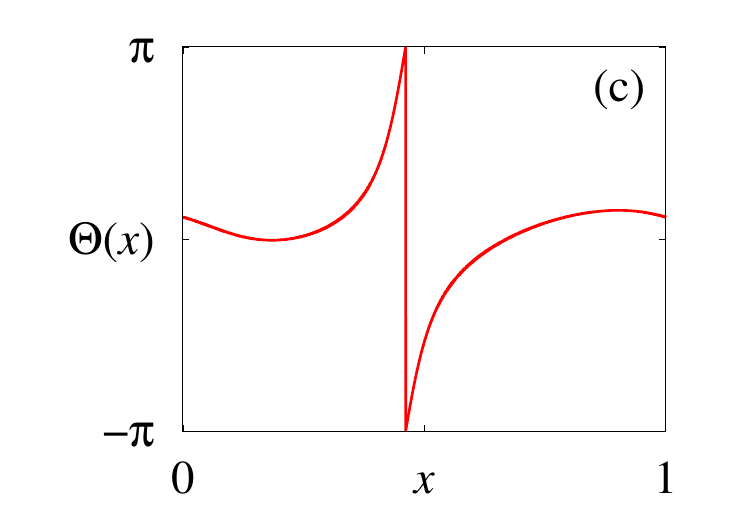}
\end{center}
\caption{The drift speed~$s$ of the solution of Eq.~(\ref{Eq:SC}) versus the parameter~$B$ for $K_3 = 0.1$. HB indicates the position of Hopf bifurcation.
Panels (a)--(e) show the arguments $\Theta(x) = \mathrm{arg}\:a(x)$
of solutions corresponding to the black dots in the main diagram.}
\label{Fig:ScanB_K0_1}
\end{figure}

\begin{figure}
\begin{center}
\includegraphics[width=0.64\textwidth]{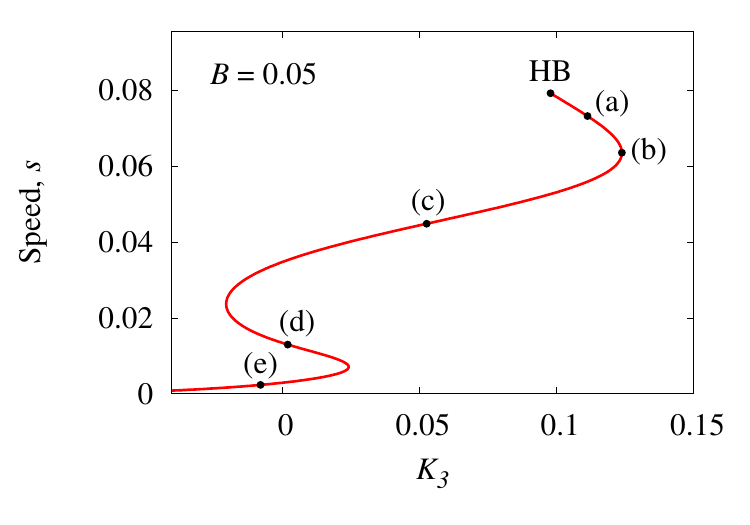}
\begin{minipage}[b]{0.32\textwidth}
\includegraphics[width=1.0\textwidth]{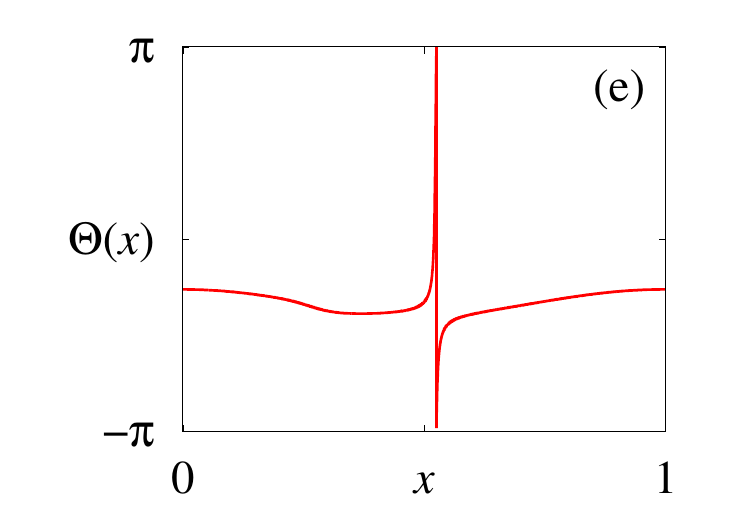}\\[2mm]
\includegraphics[width=1.0\textwidth]{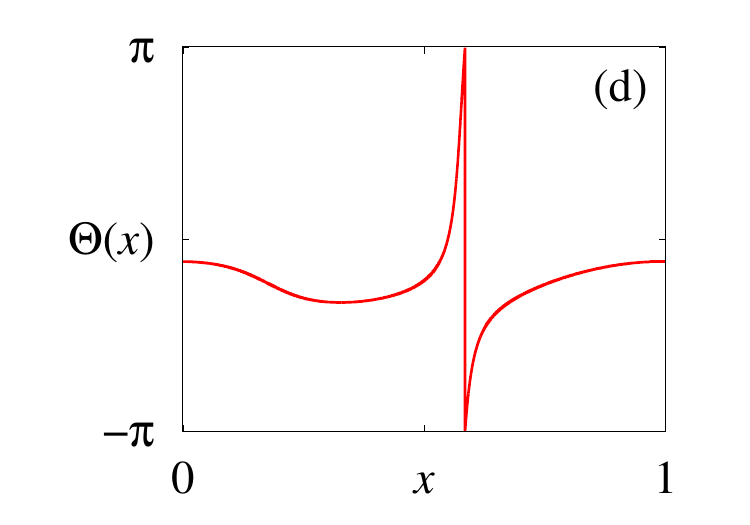}
\end{minipage}
\\[2mm]
\includegraphics[width=0.32\textwidth]{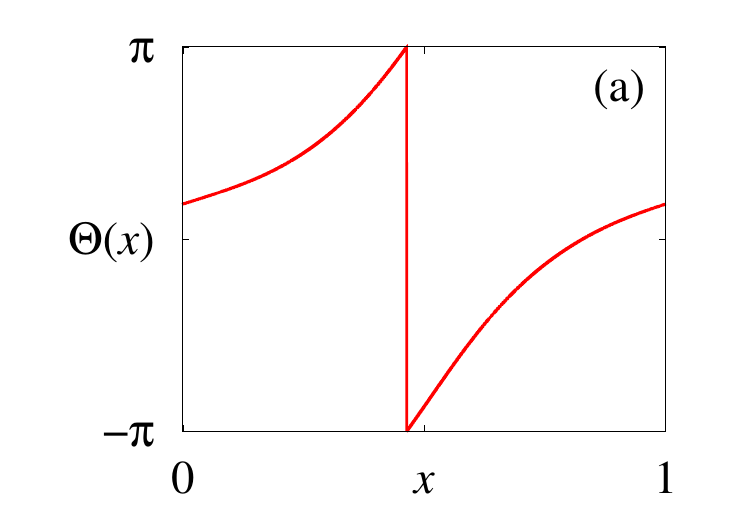}
\includegraphics[width=0.32\textwidth]{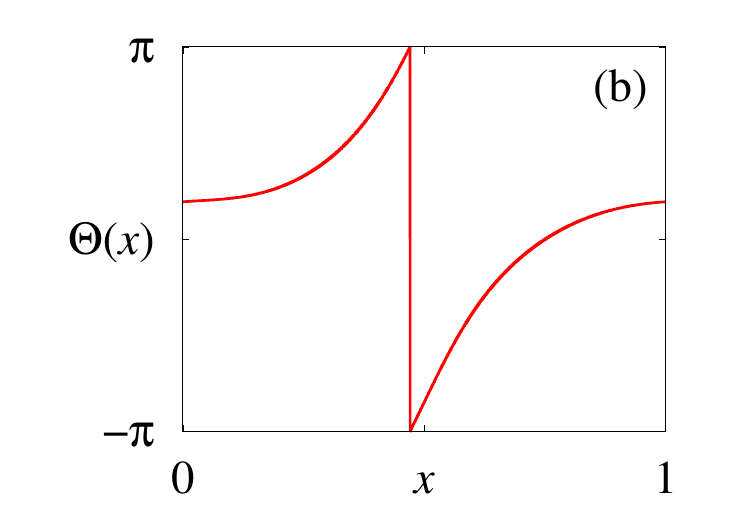}
\includegraphics[width=0.32\textwidth]{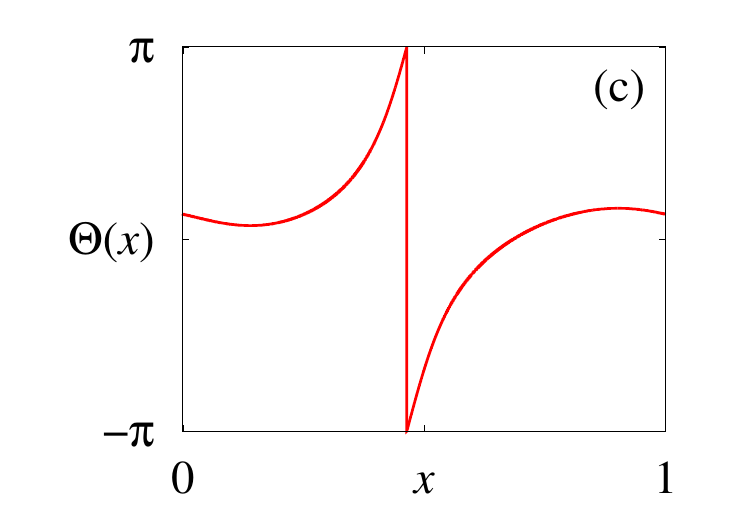}
\end{center}
\caption{The drift speed $s$ of the solution of Eq.~(\ref{Eq:SC}) versus the parameter $K_3$ for $B = 0.05$. HB indicates the position of Hopf bifurcation.
Panels (a)--(e) show the arguments $\Theta(x) = \mathrm{arg}\:a(x)$
of solutions corresponding to the black dots in the main diagram.}
\label{Fig:ScanK_B0_05}
\end{figure}

\begin{figure}
\begin{center}
\includegraphics[width=0.64\textwidth]{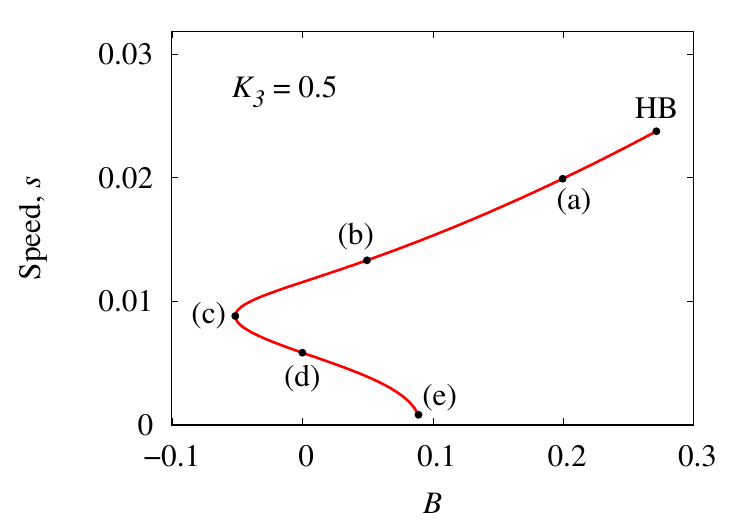}
\begin{minipage}[b]{0.32\textwidth}
\includegraphics[width=1.0\textwidth]{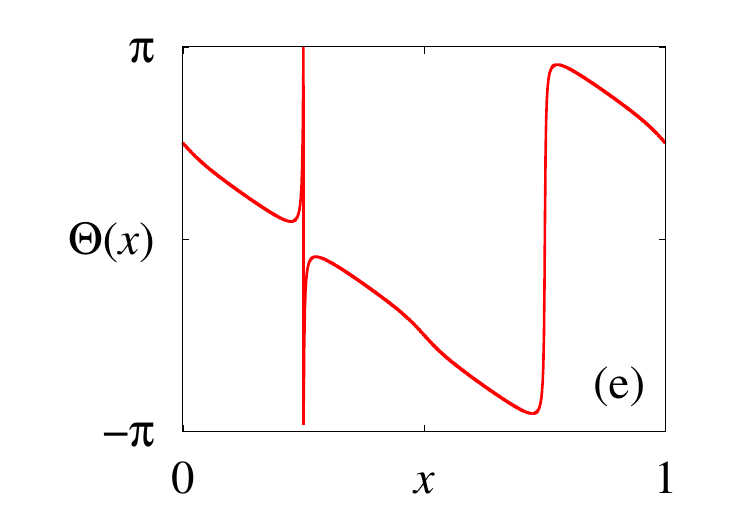}\\[2mm]
\includegraphics[width=1.0\textwidth]{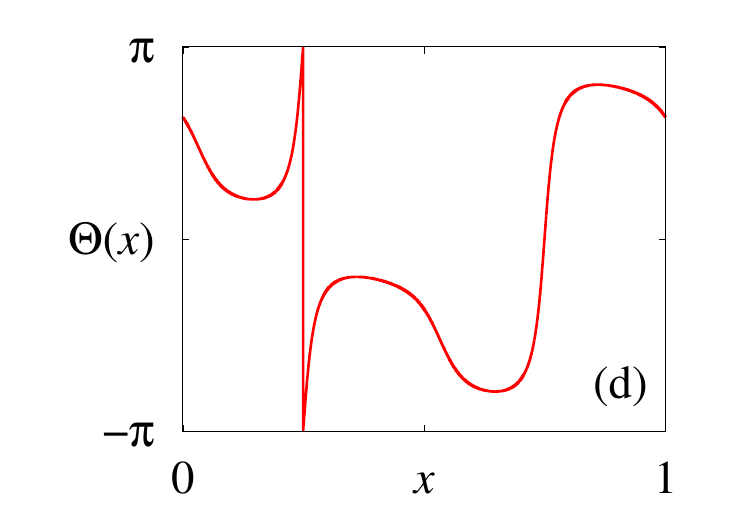}
\end{minipage}
\\[2mm]
\includegraphics[width=0.32\textwidth]{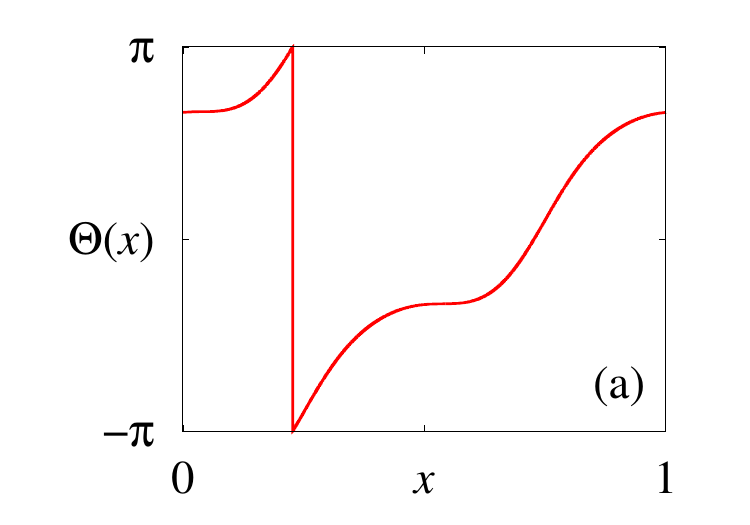}
\includegraphics[width=0.32\textwidth]{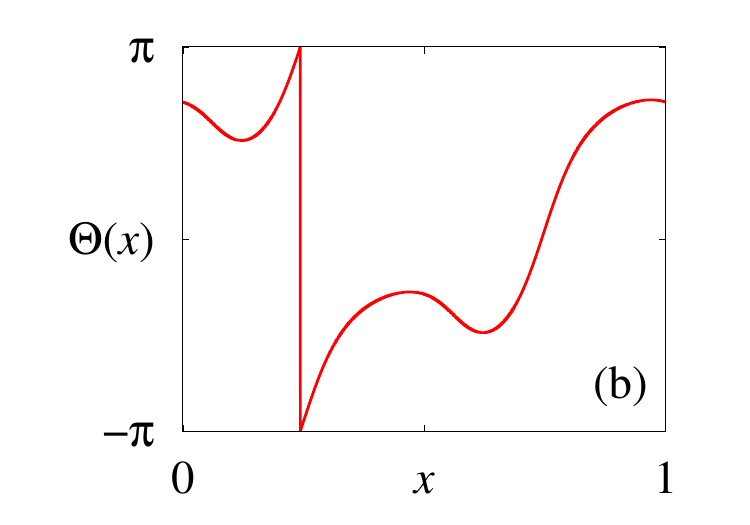}
\includegraphics[width=0.32\textwidth]{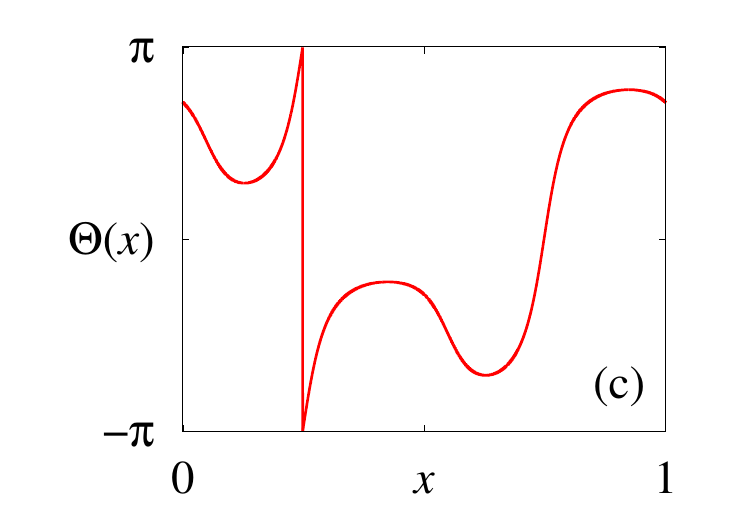}
\end{center}
\caption{The drift speed $s$ of the solution of Eq.~(\ref{Eq:SC}) versus the parameter $B$ for $K_3 = 0.5$. HB indicates the position of Hopf bifurcation.
Panels (a)--(e) show the arguments $\Theta(x) = \mathrm{arg}\:a(x)$
of solutions corresponding to the black dots in the main diagram.}
\label{Fig:ScanB_K0_5}
\end{figure}

\begin{figure}
\begin{center}
\includegraphics[width=0.64\textwidth]{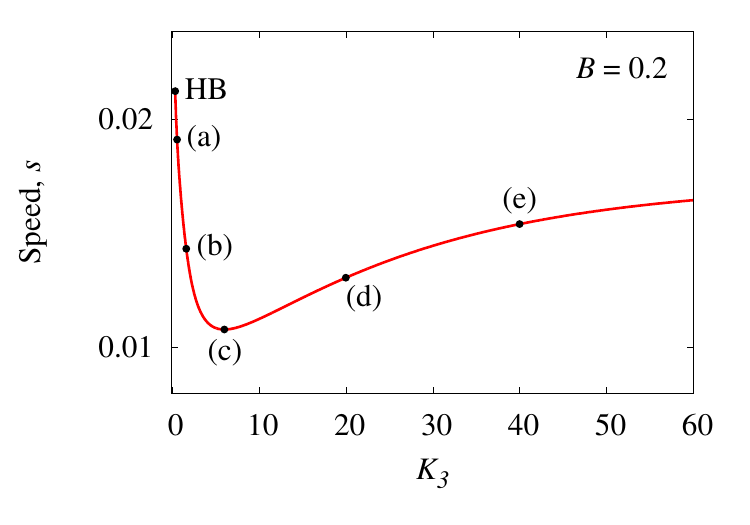}
\begin{minipage}[b]{0.32\textwidth}
\includegraphics[width=1.0\textwidth]{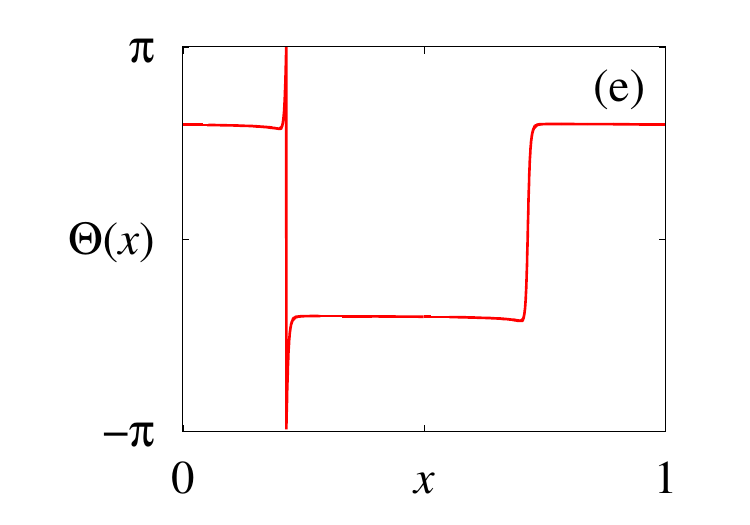}\\[2mm]
\includegraphics[width=1.0\textwidth]{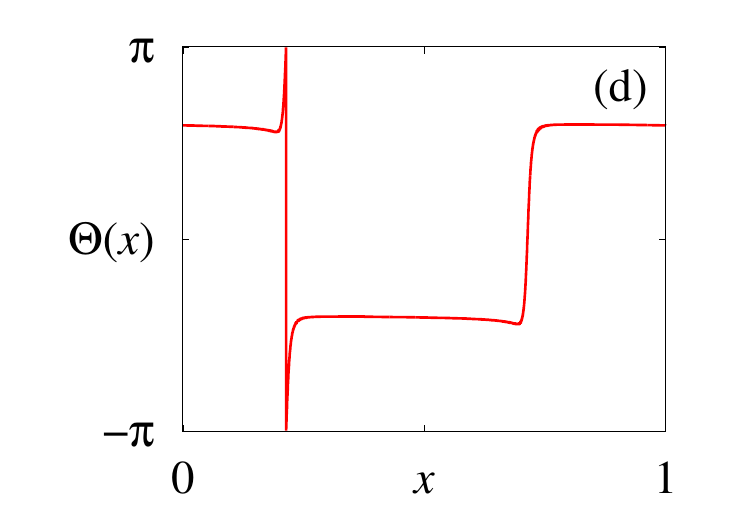}
\end{minipage}
\\[2mm]
\includegraphics[width=0.32\textwidth]{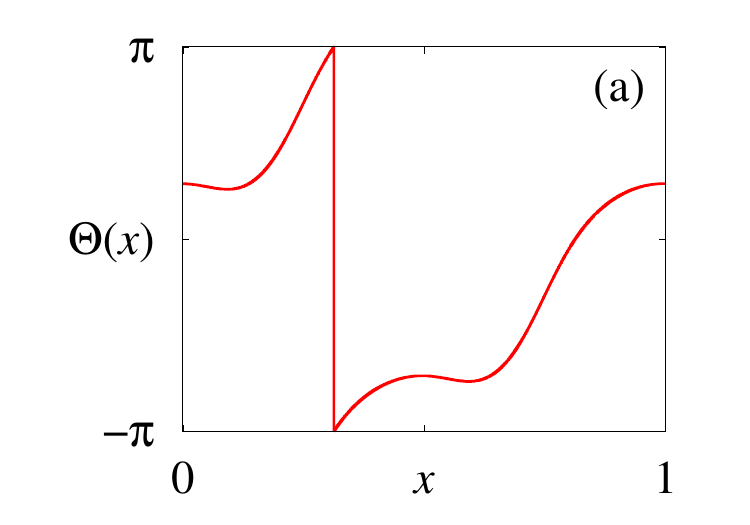}
\includegraphics[width=0.32\textwidth]{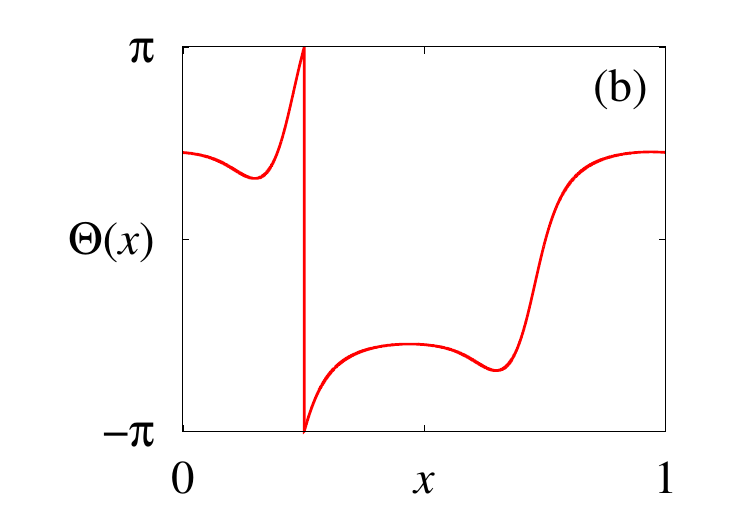}
\includegraphics[width=0.32\textwidth]{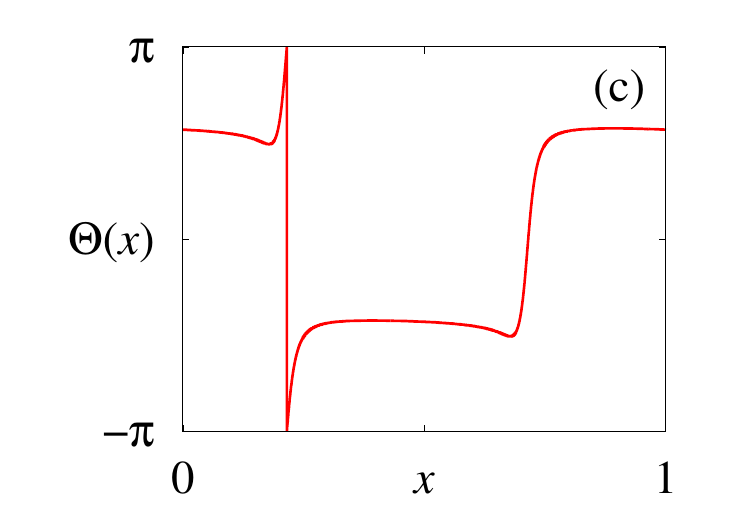}
\end{center}
\caption{The drift speed $s$ of the solution of Eq.~(\ref{Eq:SC}) versus the parameter $K_3$ for $B = 0.2$. HB indicates the position of Hopf bifurcation.
Panels (a)--(e) show the arguments $\Theta(x) = \mathrm{arg}\:a(x)$
of solutions corresponding to the black dots in the main diagram.}
\label{Fig:ScanK_B0_2}
\end{figure}

\noindent
The approach outlined in the previous section now allows to continue periodic solutions that emanate from a Hopf bifurcation of a splay solution.
The linear stability analysis of the $1$-twisted state (cf.~Section~\ref{sec:bifurcation}) indicates the location of the Hopf bifurcations; cf.~Fig.~\ref{fig:Eigenvalues} for the stability diagram of the $1$-twisted state
in the $(B,K_3)$-plane for fixed $A = 0.9$, $K_2 = 1$, $\alpha_2 = \pi/2 - 0.1$, $\alpha_3 = 0$.
Where the Hopf bifurcation is supercritical new types of stable time-dependent synchronization patterns can emerge.
We tested this theoretical prediction in numerical simulations for finite networks~\eqref{eq:finite_system} of $N = 512$
phase oscillators and found that the new solutions
take the form of spatially modulated traveling waves.
After identifying such solutions, we used the self-consistency equation~\eqref{System:w} with the pinning condition~\eqref{Pinning:w} to perform their
arc-length continuation.
Figs.~\ref{Fig:ScanB_K0_1}--\ref{Fig:ScanK_B0_2} show typical solution branches in terms of the drift speed~$s$ and the asymmetry parameter~$B$ or the strength of the higher-order interactions~$K_3$.
We do not compute stability along the branches explicitly but summarize stability properties of the branches indicated by direct numerical simulations of finite networks.

Note that for small values of~$K_3$ and~$B$, the supercritical Hopf bifurcation of $1$-twisted state is mediated by the eigenvalue $\lambda_1^+$ with the eigenfunction $v_1^+(x) = 1 - e^{-2 \pi i x}$, see Fig.~\ref{fig:Eigenvalues}(b).
For example, Figs.~\ref{Fig:ScanB_K0_1} and~\ref{Fig:ScanK_B0_05} show the solution branches for $K_3 = 0.1$ and $B = 0.05$, respectively. 
(These values correspond to the bottom horizontal and left vertical lines in Fig.~\ref{fig:Eigenvalues}.)

The branch for $K_3 = 0.1$,
as shown in Fig.~\ref{Fig:ScanB_K0_1}, 
starts at the point $(B,s)\approx(0.05,0.079)$ in which the drift speed~$s$ coincides with the value $\Imag(\lambda_1^+) / (2\pi)$ predicted by Remark~\ref{rem:Period}.
The branch consists of two parts with different slopes that meet at the fold point~(b).
Numerical simulations of the finite system~\eqref{eq:finite_system} suggest that the upper part is stable whereas the lower part is unstable.
Moving along the solution branch, we see that as the drift speed~$s$ decreases, the original straight profile~$\Theta(x)$ bends down and up in its left and right sections.
Moreover, it seems likely that the solution branch extends asymptotically to $B\to\infty$.
In this limit, the speed~$s$ vanishes and the phase profile~$\Theta(x)$ converges to a horizontal line with a phase-slip discontinuity at one point.

The shape of the branch for fixed $B = 0.05$, shown in Fig.~\ref{Fig:ScanK_B0_05}, is more complex.
It is characterized by the non-monotonic dependence
of~$s$ on~$K_3$ such that at least four different slope parts separated by three fold points can be found in the corresponding diagram.
Numerical simulations of the finite system~\eqref{eq:finite_system} suggest that the negative slope parts of the branch are stable whereas those with positive slope are unstable.
Moreover, for large negative values of~$K_3$ we observe similar asymptotic behavior as in the $B\to\infty$ limit in Fig.~\ref{Fig:ScanB_K0_1}.

Now we describe two examples of Hopf bifurcation for larger values of~$K_3$ and~$B$, when this bifurcation is mediated by the eigenvalue $\lambda_2^+$ with the eigenfunction $v_2^+(x) = 1 - e^{-4 \pi i x}$.
Note that in this case the drift speed
is equal to $\Imag(\lambda_2^+) / (4\pi)$
as expected for the linear approximation; cf.~Remark~\ref{rem:Period}.
Fig.~\ref{Fig:ScanB_K0_5} shows the solution branch for $K_3 = 0.5$; see also the top horizontal line in Fig.~\ref{fig:Eigenvalues}.
The branch starts at the point $(B,s)\approx(0.27,0.024)$, folds at the point~(c) and numerical continuation terminates at the point~(e).
Using numerical simulations of the finite system~\eqref{eq:finite_system}, we find that the lower part of the branch is unstable, whereas the upper part is only stable from the right-most point
to some point between~(b) and~(c).
This indicates that before approaching the fold point~(c), the solution is destabilized by some dynamical bifurcation, most likely a secondary Hopf bifurcation, although we have not checked this hypothesis rigorously.
At~(e) the numerical computation of the branch terminates; this is due to the fact that the coefficients~$b$ and~$e^{i\theta}$ of the M{\"o}bius transformation representing the Poincar{\'e} map of Eq.~\eqref{Eq:Riccati} satisfy the limiting relation $|b| - |\sin(\theta/2)| \to 0$.
Thus, the determinant of the self-consistency system~\eqref{System:w}, \eqref{Pinning:w} tends to infinity and numerical continuation of the solution becomes impossible.
The nature of the singularity close to point~(e) is unclear an we briefly discuss this in Section~\ref{sec:Discussion} below.

To complete the description of the solution branch shown in Fig.~\ref{Fig:ScanB_K0_5}, we note that in this case all solutions of the self consistency system~\eqref{System:w}, \eqref{Pinning:w} satisfy the identities $\hat{w}_1 = \hat{w}_4 = \hat{w}_5 = 0$. Therefore, we have
\[
w(x) = \hat{w}_2 \psi_2(x) + \hat{w}_3 \psi_3(x)
= \hat{w}_2 \cos(2\pi x) + \hat{w}_3 \sin(2\pi x).
\]
On the other hand, the profiles~$a(x)$ determined by~\eqref{Operator:U} satisfy a symmetry relation
$
a(x+1/2) = - a(x).
$
With $\Theta(x) = \mathrm{arg}\:a(x)$ this means that 
\[
\Theta(x+1/2) = \Theta(x) + \pi\:\:\mbox{mod}\:\:2\pi.
\]
The latter relation is clearly seen
in the panels (a)--(e) of Fig.~\ref{Fig:ScanB_K0_5}.

The last example of this section is the solution branch of~\eqref{System:w}, \eqref{Pinning:w} for $B = 0.2$ as shown in Fig.~\ref{Fig:ScanK_B0_2}.
It shows the nonmonotonic dependence of the drift speed~$s$ on the parameter~$K_3$.
Starting from the Hopf bifurcation point, the speed~$s$ decreases up to point~(c) and then increases for larger values of~$K_3$.
Although there are no fold points on the branch,
numerical simulations of the finite system~\eqref{eq:finite_system} show that the corresponding solutions become unstable for large values of~$K_3$.
This is another sign of the existence of secondary bifurcations in the system. 
Remarkably, the solution branch in Fig.~\ref{Fig:ScanK_B0_2} seems to have an asymptotic limit for $K_3\to\infty$.
In this limit, the speed~$s$ tends to some nonzero value and the phase profile~$\Theta(x)$ approaches a two-cluster anti-phase state (cf.~Panel~(e)) with a phase-slip discontinuity at one point.

\section{Discussion}
\label{sec:Discussion}
\noindent
In this paper, we performed a detailed analysis of Hopf bifurcations of twisted states in a ring network of phase oscillators with nonlocal higher-order interaction.
We were able to conduct not only the linear stability analysis of twisted states, but also to describe the global properties of new periodic solutions emanating from the Hopf bifurcation.
For the numerical continuation, we focused on the strength of the nonpairwise interactions~$K_3$---in line with previous work for $\alpha_2=0$~\cite{Bick2023a}---as well as the asymmetry parameter~$B$ of the coupling kernel that facilitates bifurcations to traveling solutions (see~\cite{Bick2014a,Omelchenko2019}).

\begin{figure}
\begin{center}
\includegraphics[height=0.24\textwidth]{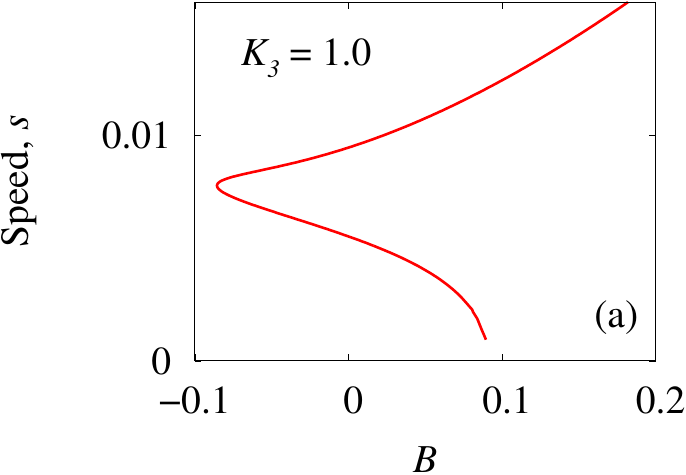}\hspace{2mm}
\includegraphics[height=0.24\textwidth]{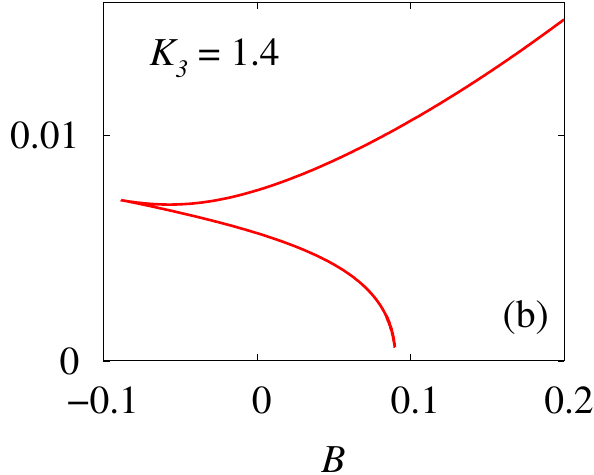}\hspace{2mm}
\includegraphics[height=0.24\textwidth]{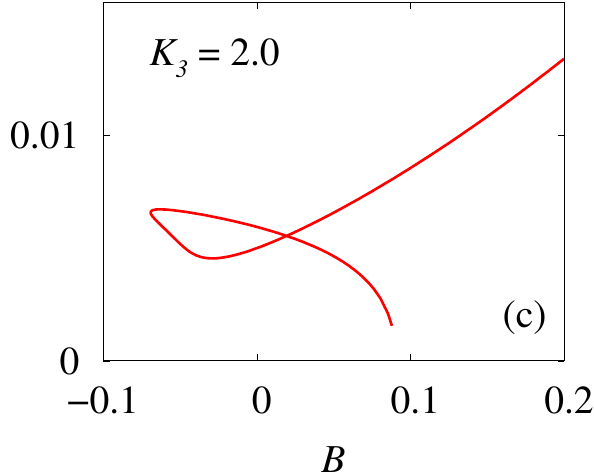}
\end{center}
\caption{
\label{fig:Cusp}
The dependence of the drift speed $s$ of traveling waves versus the parameter $B$ for (a) $K_3 = 1$, (b) $K_3 = 1.4$, and (c) $K_3 = 2$. Other parameters:
$\omega = 0$, $A = 0.9$, $K_2 = 1$, $\alpha_2 = \frac{\pi}{2} - 0.1$, and $\alpha_3 = 0$.
}
\end{figure}

While we restricted the bifurcation analysis to compute one-parameter bifurcation diagrams in Section~\ref{sec:Continuation}, computing these sheds light on potential codimension two bifurcations.
Fig.~\ref{fig:Cusp} shows one parameter bifurcation diagrams in the asymmetry parameter~$B$ for different strength~$K_3$ of the nonpairwise interactions (cf.~also Fig.~\ref{Fig:ScanB_K0_5}). For $K_3\approx1.4$ it appears that there is a cusp point where the fold point ``turns over''.
Inspecting the phase profile at the fold point for varying parameter~$K_3$ reveals that they are very similar to Fig.~\ref{Fig:ScanB_K0_5}(c).

Note that the phase profile at the end of the numerically computed solution branch in Fig.~\ref{Fig:ScanB_K0_5} for small but finite~$s$ resembles a $-1$-twisted (antisplay) phase configuration up to a single twist around the torus as shown in Panel~(e); 
similar phase profiles are obtained at the end of the branches in Fig.~\ref{fig:Cusp} (not shown). 
Indeed, the point in parameter space where numerical continuation terminates is close to a bifurcation of the (stationary) $-1$-twisted phase configuration at $B\approx0.09$: 
First, note that the linear stability of the $-1$-twisted phase configuration is given by~\eqref{eq:OneTwLinStab} with~$B$ replaced by~$-B$.
Now the real eigenvalues~$\lambda^\pm_\ell$ for $\ell\geq 4$ pass through zero at $B = A\tan(\alpha_2)^{-1}\approx 0.09$ for the parameter chosen independent of~$K_3$.
Together with the fact that~$s$ is close to zero where the numerical continuation terminates, one may speculate that these branches relate to this degenerate bifurcation of the $-1$-twisted (antisplay) phase configuration.
While this may be possible for networks of finitely many oscillators, in the limit of $N\to\infty$ the solutions are topologically different due to distinct winding numbers and the bifurcation likely involves the essential spectrum.
Clarifying the nature of the singularity requires further investigation and is beyond the scope of the current work.

In our analysis we focused on specific examples of higher-order interactions and a sinusoidal coupling kernel that facilitated the analysis. 
Many variations of the model are possible, such as considering other interactions in oscillator rings beyond Kuramoto-type pairwise interactions~\cite{Bolotov2019a,Smirnov2023} or rings of phase oscillators with nonidentical intrinsic frequencies~\cite{Omelchenko2014}.
Furthermore, turbulence is one of the complex dynamical behavior that is observed in the phase oscillator networks~\cite{Wolfrum2016} and whether this can be understood in terms of the bifurcation scenarios outlined here is an open question.
Of course, new phenomena related to twisted states can also be expected for more complicated higher-order interactions, such as adaptive higher-order interactions~\cite{Rajwani2023}.

The dynamics of phase oscillator networks naturally relate to the dynamics of more general, nonlinear oscillator networks through phase reduction (see, for example~\cite{Ashwin2015a,Leon2019a,Bick2023}). 
More specifically, higher-order phase interactions can arise through higher-order corrections to the phase reductions even when the coupling of the nonlinear oscillators is pairwise.
While rings of nonlinear oscillators beyond weak coupling can also be analyzed directly~\cite{Laing1998}, we anticipate that the bifurcation mechanisms uncovered here can shed light on the dynamics of nonlinear oscillator networks such as those discussed in~\cite{Zou2009,Lee2022}. 
Making this explicit in a specific model leaves interesting directions for future research.

\section*{Acknowledgements}
\noindent
TB and CB acknowledge support through a Hans--Fischer Fellowship at the Institute for Advanced Study of the Technische Universit\"at M\"unchen.
CB was also supported by the Engineering and Physical Sciences Research Council (EPSRC) through the grant EP/T013613/1.
Moreover, CB is grateful for the hospitality of OO at the Universit\"at Potsdam during multiple research visits.
The work of OO was supported by the Deutsche Forschungsgemeinschaft under Grant No.~OM 99/2-2.



\bibliographystyle{unsrt}
\def\urlprefix{}
\def\url#1{}

\bibliography{ref} 

\appendix

\section{Bifurcation Formulas}\label{sec:bifurcation_formulas}


{\allowdisplaybreaks

\noindent
The variable~$\zeta$ from formulas~\eqref{eq:curve_approx} can be calculated using the first, second and third derivative of~$\mathcal G^1$ at the bifurcation point. 
In order to state a formula for~$\zeta$ we continue to use the notation from Section~\ref{sec:bif}. Moreover, we denote $A:= D_v\mathcal G^1(0,p_0)$ for the linearization at the bifurcation point, $\bar v_\ell$~for the complex conjugate of~$v_\ell$, and $v_\ell'$~for an element in the dual space of~$X$ such that
\begin{align*}
    \langle v_\ell, v_\ell'\rangle = 1, \quad \text{and} \quad \langle v, v_\ell'\rangle = 0,
\end{align*}
for all other eigenfunctions of~$A$. 
According to~\cite[Formula I.9.11]{Kielhofer2012} we can then compute~$\zeta$ as
\begin{align*}
    \zeta &= -\langle D^3_{vvv} \mathcal G^1(0,p_0)[v_\ell, v_\ell, \bar v_\ell], v_\ell'\rangle\\
    &\quad - \langle D^2_{vv}\mathcal G^1(0,p_0)[\bar v_\ell, (2 i \kappa(0)- A)^{-1}D_{vv}^2\mathcal G^1(0,p_0)[v_\ell, v_\ell]], v_\ell'\rangle\\
    &\quad + 2\langle D^2_{vv} \mathcal G^1(0,p_0)[v_\ell, A^{-1}D^2_{vv}\mathcal G^1(0,p_0)[v_\ell, \bar v_\ell]], v_\ell'\rangle.
\end{align*}
If $\ell = 1$, evaluating this formula explicitly yields
\begin{align*}
	\zeta &= \Bigg\{ \frac{1}{2} K_2^2 \Big[\sin \left(2 \alpha _2\right) \left(16 A^2+59 i A B-80 A+11 B^2-48 i B+64\right)\\
	&\qquad +\cos \left(2 \alpha _2\right) \left(16 i A^2+A (5 B-16 i)+(28-33 i B) B\right)\\
	&\qquad +48 i A^2+27 A B-112 i A+65 i B^2-28 B+128 i\Big]\\
   &\quad -2 i K_3 K_2 \Big[2 \cos(\alpha_2)
   \Big[(4 \cos(\alpha_3) (4 A^3+A^2 (-3-i B)+A (-2 B+i)^2\\
   &\qquad-i B^3-3 B^2+5 i B+4) -i \sin(\alpha_3) (6 A^3+A^2 (-12+5 i B)\\
   &\qquad+A \left(5 B^2+4 i B+4\right)+6 i B^3+4 B^2-20 i B-16)\Big]\\
   &\qquad +\sin(\alpha_2) \Big[\cos(\alpha_3) \left(40 A^2 (B+i)+A (43 B+24 i)+40 B^3+51 i B^2-8 B-32 i\right)\\
   &\qquad+\sin(\alpha_3) (20 A^3+A^2 (-72+19 i B)+A \left(22 B^2-13 i B+40\right)\\
   &\qquad+25 i B^3-19 B^2+56 i B+32)\Big]\Big]\\   
   &\quad +4 K_3^2 (B-i A) \Big[-10 A^3+9 i A^2 B+8 A^2+\sin \left(2 \alpha _3\right) (6 i A^3+5 A^2 B\\
   &\qquad+A \left(5 i B^2+6 B-8 i\right)+4 B \left(B^2-2 i
   B-2\right))\\
   &\qquad+\cos(2 \alpha_3) \left(2 A^3+A^2 (-8-5 i B)+A \left(B^2-6 i B-8\right)-4 i B \left(B^2-4 i B-2\right)\right)\\
   &\qquad-9 A B^2+10 i A B+8 i B^3+12 B^2-16 i B\Big]\Bigg\}\\
   &\quad \cdot \Bigg\{ 16 \Big[K_2 \left(\cos \left(\alpha _2\right)
   (-B-i A)-\sin \left(\alpha _2\right) (A+3 i B-4)\right)\\
   &\qquad+2 K_3 (A+i B) \left(i \cos \left(\alpha _3\right)+\sin \left(\alpha _3\right) (A-i B-1)\right)\Big]\Bigg\}^{-1}.
\end{align*}
Further, if $\ell = 2$, we get
\begin{align*}
    \zeta & = \Bigg\{-\frac{1}{2} K_2^2 \Big[i \sin \left(2 \alpha _2\right) \left(A^2-4 i A B+3 B^2\right)\\
    &\qquad+\cos \left(2 \alpha _2\right) \left(A^2-2 i A B-3 B^2\right)+3 A^2+4 i A B+9 B^2\Big]\\
    &\quad +K_3 K_2 \Big[\cos \left(\alpha _2\right) \left(i \sin \left(\alpha _3\right) \left(2 A^2+9 i A B-B^2\right)+\cos \left(\alpha _3\right) \left(8 A^2+11 i A B-9 B^2\right)\right)\\
    &\qquad+\sin \left(\alpha _2\right) (B+i A) \left(11 A \cos \left(\alpha _3\right)+\sin \left(\alpha _3\right) (-10 B+7 i A)\right)\Big]\\
    &\quad +K_3^2 \Big[\cos \left(2 \alpha _3\right) \left(3 A^2+10 i A B-3 B^2\right)-9 A^2\\
    &\qquad+i \sin \left(2 \alpha _3\right) (3 A+i B)^2+8 i A B-9 B^2\Big]\Bigg\}\\
    &\quad \cdot\Big\{2 K_2 \left(\sin \left(\alpha _2\right) (B+i A)+\cos \left(\alpha _2\right) (A+i B)\right)+4 i K_3 \left(A \sin \left(\alpha _3\right)+B \cos \left(\alpha _3\right)\right)\Big\}^{-1}.
\end{align*}
Finally, when $\ell = 3$ we have
\begin{align*}
    \zeta &= \Big\{ B \left(\cos \left(\alpha _2\right)-\cos \left(\alpha _3\right)\right) \Big[2 K_2 \left(B \cos \left(\alpha _2\right)-A \sin \left(\alpha _2\right)\right)\\
    &\quad +K_3 \left(\sin \left(\alpha _3\right) \left(-A^2-6 i A B+B^2\right)+i \cos \left(\alpha _3\right) \left(3 A^2+2 i A B-3 B^2\right)\right)\Big]\Big\}\\
    &\quad \cdot \Big\{8 \left(A \cos \left(\alpha _2\right)+B \sin \left(\alpha _2\right)\right)\Big\}^{-1}.
\end{align*}
These values of~$\zeta$ can be used to determine whether a bifurcation is sub- or supercritical.
}

\section{Stability analysis of $1$-twisted states on the Ott--Antonsen Manifold}
\label{Appendix:OA}
\noindent
Each $1$-twisted state in Eq.~(\ref{eq:main_system})
corresponds to a solution
\begin{equation}z(x,t) = e^{i ( 2\pi x + \Omega t)}
\label{Ansatz:Twisted}
\end{equation}
of Eq.~\eqref{Eq:OA}.
More specifically, if we insert~\eqref{Ansatz:Twisted} into Eq.~\eqref{Eq:OA}, we obtain
\begin{equation}
i \Omega = \fr{K_2}{2} e^{i \alpha_2} \hat{g}_1 - \fr{K_2}{2} e^{-i \alpha_2} \hat{g}_{-1}
+ \fr{K_3}{2} e^{i \alpha_3} \hat{g}_2 \hat{g}_{-1} - \fr{K_3}{2} e^{-i \alpha_3} \hat{g}_{-2} \hat{g}_1,
\label{Formua:Omega}
\end{equation}
where
\[
\hat{g}_k = \int_0^1 G(x) e^{-2\pi i k x} \d x,\quad k\in\mathbb{Z},
\]
are complex Fourier coefficients of the function~$G(x)$.
Taking into account that for any real function~$G(x)$ we have
\[
\overline{\hat{g}_k} = \hat{g}_{-k},
\]
the relation~\eqref{Formua:Omega} can be rewritten as
\begin{equation}
\Omega = K_2 \Imag\left( e^{i \alpha_2} \hat{g}_1 \right) + K_3 \Imag\left( e^{i \alpha_3} \hat{g}_2 \hat{g}_{-1} \right).
\label{Formua:Omega_}
\end{equation}
The latter is a formula expressing the relationship
between the frequency of $1$-twisted states
and the coupling parameters between oscillators.

To perform a linear stability analysis of the solution~\eqref{Ansatz:Twisted}, it is convenient to transform Eq.~\eqref{Eq:OA} into a co-rotating frame
\[
z(x,t)\mapsto u(x,t),\quad\mbox{where}\quad
z(x,t) = u(x,t) e^{i ( 2\pi x + \Omega t)}.
\]
After this transformation, we obtain
\begin{equation}
    \begin{split}
\pf{}{t} u(x,t) &= - i \Omega u + \fr{K_2}{2} e^{i \alpha_2} e^{-2\pi i x} \mathcal{G}\left( e^{2\pi i x} u \right)\\ 
&\qquad -\fr{K_2}{2} e^{-i \alpha_2} u^2 e^{2\pi i x} \mathcal{G}\left( e^{-2\pi i x} \overline{u} \right)\\
&\qquad+ \fr{K_3}{2} e^{i \alpha_3} e^{-2\pi i x} \mathcal{G}\left( e^{4\pi i x} u^2 \right)\:\mathcal{G}\left( e^{-2\pi i x} \overline{u} \right)\\
&\qquad -\fr{K_3}{2} e^{-i \alpha_3} u^2 e^{2\pi i x} \mathcal{G}\left( e^{-4\pi i x} \overline{u}^2 \right)\:\mathcal{G}\left( e^{2\pi i x} u \right).
\end{split}\label{Eq:OA_}
\end{equation}
Inserting the ansatz
\[
u(x,t) = 1 + v(x,t)
\]
into Eq.~\eqref{Eq:OA_} and linearizing the result
with respect to the small perturbation~$v$, we obtain
\begin{equation}
    \begin{split}
\pf{}{t} v(x,t) &= - \eta v + \fr{K_2}{2} e^{i \alpha_2} e^{-2\pi i x} \mathcal{G}\left( e^{2\pi i x} v \right) - \fr{K_2}{2} e^{-i \alpha_2} e^{2\pi i x} \mathcal{G}\left( e^{-2\pi i x} \overline{v} \right) \\
&\qquad+ K_3 e^{i \alpha_3} \hat{g}_{-1} e^{-4\pi i x} \mathcal{G}\left( e^{4\pi i x} v \right)
+ \fr{K_3}{2} e^{i \alpha_3} \hat{g}_2 e^{2\pi i x} \mathcal{G}\left( e^{-2\pi i x} \overline{v} \right) \\
&\qquad- K_3 e^{-i \alpha_3} \hat{g}_1 e^{4\pi i x} \mathcal{G}\left( e^{-4\pi i x} \overline{v} \right)
- \fr{K_3}{2} e^{-i \alpha_3} \hat{g}_{-2} e^{-2\pi i x} \mathcal{G}\left( e^{2\pi i x} v \right).
    \end{split}
    \label{Eq:OA:Lin}
\end{equation}
where
\[
\eta = i \Omega + K_2 e^{-i \alpha_2} \hat{g}_{-1} + K_3 e^{-i \alpha_3} \hat{g}_{-2} \hat{g}_1.
\]

To investigate the decay of different spatial Fourier modes, we insert the ansatz
\[
v(x,t) = v_+ e^{2\pi i k x} e^{\lambda t} + \overline{v}_- e^{-2 \pi i k x} e^{\overline{\lambda} t}
\quad\mbox{with}\quad v_+,v_-\in\mathbb{C}\:\:\mbox{and}\:\: k\in\mathbb{Z}
\]
into Eq.~\eqref{Eq:OA:Lin} and equate separately the terms at~$e^{\lambda t}$ and $e^{\overline{\lambda} t}$.
Thus, we obtain a spectral problem
\[
\lambda \left(
\begin{array}{c}
v_+ \\[2mm]
v_-
\end{array}
\right) =
\mathbf{B}(k)
\left(
\begin{array}{c}
v_+ \\[2mm]
v_-
\end{array}
\right),
\]
where
\[
\mathbf{B}(k) = \left(
\begin{array}{cc}
-\eta + p(k) & -q(k) \\[2mm]
-p(k) & -\overline{\eta} + q(k)
\end{array}
\right)
\]
and
\begin{align*}
p(k) &= \fr{K_2}{2} e^{i \alpha_2} \hat{g}_{k+1} + K_3 e^{i \alpha_3} \hat{g}_{-1} \hat{g}_{k+2}
- \fr{K_3}{2} e^{-i \alpha_3} \hat{g}_{-2} \hat{g}_{k+1},\\[2mm]
q(k) &= \fr{K_2}{2} e^{-i \alpha_2} \hat{g}_{k-1} + K_3 e^{-i \alpha_3} \hat{g}_1 \hat{g}_{k-2}
- \fr{K_3}{2} e^{i \alpha_3} \hat{g}_2 \hat{g}_{k-1}.
\end{align*}
The corresponding characteristic equation reads
$\mathrm{det}( \lambda \mathbf{I} - \mathbf{B}(k) ) = 0$.
It can be solved explicitly, which yields
\[
\lambda_\pm(k) = \fr{1}{2} \left( \mathrm{tr}\:\mathbf{B}(k) \pm \sqrt{ \left( \mathrm{tr}\:\mathbf{B}(k) \right)^2 - 4 \mathrm{det}\:\mathbf{B}(k) } \right).
\]

Note that for the coupling function~\eqref{Formula:G} only three leading coefficients are non-zero
\[
\hat{g}_0 = 1,\quad \hat{g}_1 = \fr{A - B i}{2},\quad \hat{g}_{-1} = \fr{A + B i}{2},
\]
while $\hat{g}_k = 0$ for $|k|\ge 2$. These relations imply
\[
\eta = i K_2 \Imag\left( e^{i \alpha_2} \hat{g}_1 \right) + K_2 e^{-i \alpha_2} \hat{g}_{-1} = K_2 \Real\left( e^{i \alpha_2} \hat{g}_1 \right)
\]
and
\[
\lambda_\pm(k) = - \eta = - \fr{K_2}{2} ( A \cos(\alpha_2) + B \sin(\alpha_2) )
\quad\mbox{for all}\quad |k| \ge 4.
\]
As for the other eigenvalues, they coincide with the expressions obtained in Section~\ref{sec:Lin}.
Moreover, a pair of zero eigenvalues is also determined by the formulas of~$\lambda_\pm(k)$.

\end{document}